\documentclass[final]{siamltex}

\usepackage{graphics}

\usepackage{amsmath,amssymb,amsfonts,amsthm}
\usepackage{mathrsfs}
\usepackage{wasysym}
\usepackage{algorithm}
\usepackage{algorithmic}

\usepackage{enumitem}
\usepackage{bm,bbm}
\usepackage{cite}

\usepackage{multirow,booktabs,comment}
\usepackage{epstopdf}

\usepackage{todonotes}
\usepackage[normalem]{ulem}

\usepackage{xcolor}
\usepackage{soul}

\usepackage[colorlinks=true,urlcolor=blue,citecolor=red,
anchorcolor=black,linkcolor=blue]{hyperref}

\usepackage{subcaption}

\usepackage{lineno}
\newtheorem{thm}{Theorem}
\newtheorem{lem}[thm]{Lemma}

\theoremstyle{definition}

\newtheorem{re}[thm]{Remark}

\newcommand{\bx}{\boldsymbol{x}}
\newcommand{\by}{\boldsymbol{y}}
\newcommand{\bX}{\boldsymbol{X}}

\newcommand{\bu}{\boldsymbol{u}}
\newcommand{\bv}{\boldsymbol{v}}

\newcommand{\bI}{\boldsymbol{I}}

\newcommand{\bXheta}{\boldsymbol{\theta}}

\usepackage{lineno}
\begin{document}
	
\title{Moving sample method for solving time-dependent partial differential equations}
\author{Beining Xu\thanks{First author. The School of Information Science and Technology and
The Institute of Mathematical Sciences, ShanghaiTech
University, Shanghai 201210 (\tt xubn2024@shanghaitech.edu.cn)}
\and
Haijun Yu\thanks{State Key Laboratory of Mathematical Sciences (SKLMS) \& LSEC, Institute of Computational Mathematics and Scientific/Engineering Computing, Academy of Mathematics and Systems Science, Chinese Academy of Sciences, Beijing 100190, China; School of Mathematical Sciences, University of Chinese Academy of Sciences, Beijing 100049, China
 (\tt hyu@lsec.cc.ac.cn)}
\and
Jiayu Zhai\thanks{Corresponding author. The Institute of Mathematical Sciences,
ShanghaiTech University, Shanghai 201210 (\tt zhaijy@shanghaitech.edu.cn)}
\and
Kejun Tang\thanks{Department of Mathematics, School of Sciences, Great Bay University, Dongguan 523000, China(\tt tangkj@gbu.edu.cn)}
\and
Xiaoliang Wan\thanks{Department of Mathematics,
Center for Computation and Technology, Louisiana State
University, Baton Rouge 70803 (\tt xlwan@math.lsu.edu)}
}

\maketitle

\begin{abstract}
Solving time-dependent partial differential equations (PDEs) that exhibit sharp gradients or local singularities is computationally demanding, as traditional physics-informed neural networks (PINNs) often suffer from inefficient point allocation that wastes resources on regions already well-resolved. This paper presents an adaptive sampling framework for PINNs aimed at efficiently solving time-dependent partial differential equations with pronounced local singularities. The method employs a residual-driven strategy, where the spatial–temporal distribution of training points is iteratively updated according to the error field from the previous iteration. This targeted allocation enables the network to concentrate computational effort on regions with significant residuals, achieving higher accuracy with fewer sampling points compared to uniform sampling. Numerical experiments on representative PDE benchmarks demonstrate that the proposed approach improves solution quality.
\end{abstract}

\section{Introduction}
%For low-dimensional partial differential equations (PDEs), grid-based numerical methods, e.g. finite element method (FEM), are commonly used efficient methods. However, when dealing with high-dimensional partial differential equations, it requires mesh generation and solving large algebraic systems on discretized grids. As a more flexible mesh-free method, neural network methods (NN) were introduced to approximate the solutions of the equations, addressing the aforementioned difficulties through machine learning methods. Its advantage lies in not requiring mesh generation; instead, automatic differentiation is performed directly on the parameterized neural network function at randomly sampled collocation points. By constructing a suitable loss function, the neural network is made to approximate the solution function of the differential equation. 
For low-dimensional partial differential equations (PDEs), traditional grid-based numerical methods, such as the finite element method (FEM), finite difference method, spectral method, etc, have proven highly effective and reliable. However, these approaches become computationally prohibitive for high-dimensional PDEs due to the curse of dimensionality, which entails exponential growth in the number of grid points required for mesh generation and the solution of large algebraic systems. To circumvent these challenges, mesh-free neural network-based methods, particularly physics-informed neural networks (PINNs), have emerged as powerful alternatives. These methods parameterize the solution via a neural network and leverage automatic differentiation to enforce the governing PDE constraints (directly through strong-form residuals or indirectly via variational/weak forms) at randomly sampled collocation points, eliminating the need for mesh generation. By minimizing a suitably constructed loss function that incorporates the PDE residual, initial conditions, and boundary conditions, the trained neural network yields an accurate approximation of the PDE solution. Typical examples include Physics-Informed Neural Networks (PINNs) \cite{raissi2019physics}, the Deep Ritz method \cite{e2018deepritz}, Weak Adversarial Networks (WANs) \cite{zang2022wan}, and so on. 
%In other words, the physical information of the differential equation is encoded into the neural network through the loss function, unlike numerical methods like FEM which solve by establishing relationships between function values at grid points through numerical differentiation.

%Furthermore, some studies indicate that, compared to grid-based methods like FEM, it requires relatively few collocation points to learn the PDE 
Although neural networks can achieve reasonable PDE approximations using fewer collocation points than traditional grid-based methods such as FEM (see e.g. \cite{e2018deepritz, zhai2022deep}), 
%This advantage is particularly evident for high-dimensional problems. However, 
a critical limitation arises when the solution contains singularities or regions of large gradients, which are typically confined (close) to some lower-dimensional manifolds within the domain. Because these singular regions occupy a much smaller volume compared to the surrounding smooth areas, insufficient collocation points are naturally sampled there under uniform or quasi-random distributions. As a result, the optimization process ---driven by the average residual over all points--- receives only weak error signals from the singular regions and predominantly focuses computational effort on reducing errors in the larger smooth regions, leading to poor resolution of sharp features and inaccurate overall solutions. This phenomenon is even more serious in high dimensional problems, where neural network methods are more efficient, particularly when singularities and their locations are unknown. Therefore, it is essential to rationally allocate these limited collocation points by employing adaptive sampling strategies that deliberately concentrate points in regions of low regularity or high residual. %On the other hand, if singularity exists in unknown regions (usually in a lower dimensional manifold), the cumulative error in these regions is relatively insignificant compared with that of the complement. Then the learning process will allocate most of the computational resources to the smooth regions while neglecting singular regions, leading to erroneous results. So we need to utilize these relatively few points rationally. %From a computational perspective, these points need to be sampled where the PDE has prominent features, such as low regularity. 

%In contrast, the more mature adaptive finite element method or moving mesh method yields satisfactory results. Adaptive numerical methods aim to generate more grid points in regions with larger errors (typically where the solution function varies sharply) when solving PDEs with FEM. This allows the iterative solution process to simultaneously adaptively adjust the mesh generation strategy, enabling the solution of more complex PDEs with a smaller number of grid points. Another typical adaptive method is the moving mesh method, used primarily for solving time-dependent PDEs. Its main idea is to control the movement of grid points by constructing control equations. While maintaining the number of grid nodes and their connectivity, it couples the mesh within the solution domain with the solution function, allowing grid nodes to change with time and continuously move to areas where the solution function has low regularity, requiring more computational resources for approximation, thus solving PDEs efficiently. Due to their significant improvement in computational efficiency, they quickly became a mainstream research direction in FEM.
Adaptive techniques have been well developed for traditional numerical methods, such as adaptive finite element methods ($h$-version) and moving mesh methods ($r$-version). Adaptive finite element methods dynamically coarsen or refine the mesh based on an indicator or estimator of local errors, thereby effectively capturing complex solution features using significantly fewer grid points overall. Another prominent adaptive strategy is the moving mesh method, which is especially effective for time-dependent PDEs. This approach fixes the number of grid nodes and their connectivity but redistributes them over time by solving auxiliary mesh-control equations that couple the grid motion to the evolving solution. As a result, nodes continuously migrate toward regions that demand greater resolution, enabling efficient computation without increasing the total degrees of freedom.

As a counterpart to adaptive grid-based methods, a series of adaptive sampling methods has been developed for neural network approaches. For example, Residual-based Adaptive Refinement (RAR) \cite{lu2021DeepXDE, wu2023rar} borrows the idea of mesh refinement based on a posteriori error estimates from adaptive FEMs. It continuously computes the total residual of the loss function in local regions during the neural network solution process and further samples, increasing the number of collocation points in regions where the residual exceeds a threshold. In \cite{gu2021Selectnet}, SelectNet compares the residual at each point to a threshold during training to decide whether to use that point in the next training phase or to determine its weight in the loss functional calculation. Failure-Informed based Adaptive Sampling (FI-PINNs) \cite{gao2023fipinn, gao2024fipinn1} uses a threshold measurement to determine where new collocation points need to be sampled. The idea behind this type of method is similar to adaptive mesh refinement in classical numerical methods. They can indeed achieve adaptive sampling, but it is difficult to quickly and directly determine the proportional number of additional points required in different local regions and allocate the increased collocation points accordingly.

Another idea for adaptive collocation point sampling is to treat the (normalized) error distribution (often the residual, sometimes including gradient information) across the entire solution domain as a probability distribution. Sampling points according to this distribution can also achieve an adaptive effect. For instance, \cite{gao2023Active} uses the Metropolis-Hastings (MH) algorithm to handle the error distribution without requiring normalization. The Deep Adaptive Sampling (DAS) method was proposed in \cite{tang2023das}, which trains a normalizing flow to generate adaptive collocation point samples. Gaussian or Gaussian mixture models are used as sampling models in \cite{hou2023enhancing, wang2025Navigating, jiao2024GAS}. For training the sampling models, \cite{tang2023adversarial, hou2023enhancing, wang2025Navigating, jiao2024GAS, han2022residual} used adversarial training. 

The aforementioned studies all design adaptive methods from a computational perspective; yet none explain the adaptive effect or prove their convergence mathematically. In contrast, \cite{tang2023adversarial} pioneered the use of optimal transport theory to demonstrate the constraints and prove the convergence of the Adversarial Adaptive Sampling (AAS) method. It frames the adaptive solving process as a probability transport process of learning collocations, where the target probability distribution corresponds to the adaptive collocation distribution in the residual sense. Using these collocation points for solving yields results whose associated residuals follow a uniform distribution, thereby achieving the adaptive effect.

For time-dependent equations, simply treating the time variable as an extra spatial dimension for extension can cause the loss of inherent dynamic information in the equation; at the same time, performing adaptive sampling directly with this treatment can lead to an excessive reduction of collocation points in certain time ranges, which is unfavorable to the overall solution of the equation. Therefore, we need to treat the spatial and temporal variables separately. There are some works following this approach. The normalized residual distribution at discrete time steps is treated as the invariant distribution of stochastic differential equations (SDEs), then these SDEs are used to generate adaptive samples accordingly~\cite{Yuxiao_Coupling_2023, Bruna2024Neural, WangHu2024NonuniformRandom, LiuZDYu2026}. However, the problem is that it requires subdividing each time interval further, allowing the SDEs to obtain samples from the invariant distribution. This imposes relatively high requirements on the design of the SDE and its solver to ensure its transient distribution can quickly evolve to the invariant distribution (see e.g. \cite{dobson2021using}). Similarly, \cite{gao2023Active} uses MCMC at each time step to generate adaptive samples.

The existing methods decouple the dynamics of the PDEs from those of the adaptive samples. In this work, we find the dynamics of adaptive collocation points over time based on the time-varying residual or gradient distribution of the neural network solution, and then use it for sampling. In fact, this is similar to the dynamics of the grid in the moving mesh method 
\cite{huang2010adaptive,Herbst1983Equidistributing,huang1994MMPDES,white1979selection,LiTZ2001MovingMesh}. We note that the work of MMPDE-net \cite{YangYDH2024MovingSampling} also employs the idea of the moving mesh method to adjust the sample of PINNs, but our approach is significantly different. The method proposed in this work can be regarded as a probabilistic version of the moving mesh method and is called moving sample method (MSM). More specifically, the velocity field of sample movements will be learned through its continuity equation, which contains the dynamical information of the residual over time. The dynamics can be easily obtained after each iteration of learning since the solution is approximated with a neural network, and the residual and the gradient are naturally available. Thus, the dynamics of the time-dependent PDE and the adaptive samples can be coupled in a more natural way. 

The remaining part of this paper is organized as follows. We recall the PINNs in Section \ref{PINNs} and derive the formulation of MSM in Section \ref{dynamics}. The algorithm implementation of MSM is introduced in Section \ref{Algorithm Implementation}, and its performance is demonstrated in Section \ref{Numerical Experiments} on several examples. We summarize the paper in Section \ref{conclusion}.

\section{Introduction to Physics-Informed Neural Networks}\label{PINNs}
We consider solving a general time-dependent differential equation problem:
\begin{equation}\label{timePDE}
\left\{
    \begin{array}{ll}
        \mathcal{L}u  = f(\bx,t) & (\bx, t)\in \Omega \times [0,T],\\
        \mathbf{\mathcal{I}}u(\bx, 0) = \bu_0(\bx) & \bx\in\Omega,\\
        \mathcal{B}u(\bx, t) = u_b(\bx, t) & (\bx, t) \in \partial \Omega \times [0, T],
    \end{array}
\right.
\end{equation}
where $u$ is in some function class $\mathcal{U}$ of functions on $\Omega \times [0,T]$, $f$ is the nonhomogeneous term, $\mathbf{\mathcal{I}}$ is the initial operator on $\Omega$ (a vector including initial values, velocities, etc, depending on the problem), $\mathcal{B}$ is the boundary operator on $\partial \Omega$, $\bu_0(\bx)$ is the initial condition, and $u_b(\bx, t)$ is the boundary condition. The differential operator $\mathcal{L}=\mathcal{L}_{(\bx,t)}$ can include both spatial and temporal differential operators, and can be dependent on spatial variables $\bx$ and time $t$ respectively.

PINNs~\cite{raissi2019physics} provide a unified framework for solving such partial differential equations by embedding physical laws directly into the training objective of neural networks. The unknown solution $u(\bx, t)$ is approximated by a parameterized deep neural network $u_{\bXheta}(\bx, t)$. The network is constructed as a composition of linear transformations and nonlinear activation functions, resulting in a nested structure:
\begin{equation}
    u_{\bXheta}(x) = W^{(L)} \, \sigma\!\Big( W^{(L-1)} \, \sigma\!\big( \cdots \sigma\!\big( W^{(1)} x + b^{(1)} \big) \cdots \big) + b^{(L-1)} \Big) + b^{(L)}.
\end{equation}
Here, $\sigma(\cdot)$ denotes the activation function, while $W^{(l)}$ and $b^{(l)}$ are the weight matrices and bias vectors of the $l$-th layer. The collection of all trainable parameters is denoted by $\bXheta=\{W^{(l)}, b^{(l)}\}_{l=1}^L$. The loss function is generally formulated as a weighted combination of the PDE residual term, the boundary condition term, and the initial condition term: 
\begin{equation}\label{loss_pinns}
    J(u_{\bXheta}) = J_{\text{PDE}}(u_{\bXheta}) + \beta_1 J_{\text{IC}}(u_{\bXheta}) + \beta_2 J_{\text{BC}}(u_{\bXheta}),
\end{equation}
where $\beta_1, \beta_2 > 0$ are the balancing weights of the loss function. For equation \eqref{timePDE}, the residual-based loss is defined as \begin{equation}
    J_{\text{PDE}}(u_{\bXheta}) = \int_0^T \int_{\Omega} |\mathcal{L}u_{\bXheta} - f|^2 d\bx dt,
\end{equation}
and the initial loss is defined as \begin{equation}
    J_{\text{IC}}(u_{\bXheta}) = \int_\Omega |\mathbf{\mathcal{I}}u_{\bXheta} - \bu_0|^2 d\bx.
\end{equation}
If the equation involves a boundary condition, a corresponding boundary loss $J_{BC}$ can be defined in a similar manner:
\begin{equation}
    J_{\text{BC}}(u_{\bXheta}) = \int_0^T\int_{\partial \Omega} |\mathcal{B}u_{\bXheta} - u_b|^2 d\bx dt.
\end{equation}

To optimize the loss function \eqref{loss_pinns}, the vanilla PINNs method numerically discretizes the integrals involved. Let $S = \{(\bx_\Omega^{(i)}, t^{(i)})\}_{i=1}^{N}$, $S_0 = \{\bx_0^{(i)}\}_{i=1}^{N_0}$ and $S_b = \{(\bx_{\partial\Omega}^{(i)}, t^{(i)})\}_{i=1}^{N_b}$ be the sets of uniformly distributed collocation points for the PDE loss $J_{\text{PDE}}$, the initial loss $J_{\text{IC}}$ and the boundary loss $J_{\text{BC}}$, respectively. Then the original functional loss can be replaced by the following empirical loss: \begin{align}
    I(u_{\bXheta}) = & \frac{1}{N} \sum_{i=1}^N |\mathcal{L}u_{\bXheta}(\bx_\Omega^{(i)}, t^{(i)}) - f(\bx_\Omega^{(i)}, t^{(i)})|^2\notag\\
    & + \beta_1 \frac{1}{N_0} \sum_{i=1}^{N_0} |\mathbf{\mathcal{I}}u_{\bXheta}(\bx_0^{(i)}, 0) - u_0(\bx_0^{(i)})|^2\notag\\
    & + \beta_2 \frac{1}{N_b} \sum_{i=1}^{N_b} |\mathcal{B}u_{\bXheta}(\bx_{\partial\Omega}^{(i)}, t^{(i)}) - u_b(\bx_{\partial\Omega}^{(i)}, t^{(i)})|^2.
\end{align}
This discretization transforms continuous objectives into tractable optimization ones, ensuring that the network solution $u_{\bXheta}$ adheres to the underlying physical principles at the sampled collocation points. In practice, gradient-based algorithms such as Adam or L-BFGS are then employed, leveraging automatic differentiation to efficiently compute the derivatives of the residual terms.

In this work, we will focus on the first term -- the PDE part, since the other terms can be similarly treated. For convenience, we denote the residual of the PDE as
\begin{equation}\label{notation_residual}
    r(\bx,t) = \mathcal{L}u(\bx,t) - f(\bx,t),
\end{equation}
for any function $u$ in the function class $\mathcal{U}_{\text{NN}}$ of neural networks.

\section{Moving Sample Method (MSM)}\label{MSM}

\subsection{The dynamics of adaptive sampling -- a probability perspective of moving-mesh methods
}\label{dynamics}

In \cite{tang2023das} and \cite{tang2023adversarial}, the necessity of adaptive sampling for PINNs is explained from the perspective of reducing statistical error through variance reduction. Thus, the collocation points $\{\bx_i\}_{i=1}^{N}$ need to be sampled accordingly to make the residual $r_t$ evenly distributed on $\Omega$, which corresponds to a probabilistic perspective of the equidistribution principle in moving-mesh methods \cite{white1979selection, huang2010adaptive, huang1994MMPDES}. As in the moving mesh methods, this sampling distribution can be chosen as $r_t^{2\gamma}$, $|u|^{\gamma}$, $|\nabla_x u|^{\gamma}$, or their mixture. 

In this work, we mainly choose unnormalized $r_t^{2}+\varepsilon$ as the sampling distribution, where $\varepsilon>0$ is a small constant. The reason of using a constant $\varepsilon>0$ is that
\begin{enumerate}
    \item[(1)] from computational perspective, as a ``preconditioner'', it makes the computation stable since its integral may be on the denominator of the formulation (see Section \ref{Subsec:MSM});
    \item[(2)] from approximation perspective, it provides a proportion of uniformly distributed samples that can help the neural network learn the global (away from the singularity region) solution well. 
\end{enumerate}
The reason for choosing the residual $r_t^{2}$ is that it can be continuously accessed throughout the learning process. Furthermore, it is usually not practical to use gradient at the beginning, since the neural network approximation can have wrong gradient information. For example, the Burgers' equation with high singularity in Section \ref{burgers}, the approximation may be almost a constant and cannot catch the correct gradient distribution at the beginning. So it is usually not efficient at the beginning. This can be overcome by using a mixture of them with an increasing penalty on the gradient part. %The residual $r_t$ is also separately computed at need first in \cite{Yuxiao_Coupling_2023}. Similarly, we don't consider here methods to update adaptivity like adversarial adaptive sampling (AAS) \cite{tang2023adversarial}, and accordingly its sampling, and consider the dynamics of the residual evolution in time variable $t$. We will further combine methods like AAS to it. 

In \cite{tang2023adversarial}, it is discussed that by changing the sampling distribution in PINNs, the adaptive effect can be achieved. So in this work, we shall employ a formulation of PINNs with respect to some dynamical measure $\mu_t$ instead of the uniform distribution to solve \eqref{timePDE}, that is to solve for all $t\in [0,T]$
\begin{equation}\label{timePDENN}
    \min_{\bXheta}\Big[ \int_\Omega |\mathcal{L}u_{\bXheta}(\bx, t) - f(\bx, t)|^2 \,\mu_t(d\bx) \approx 
    \frac{1}{N} \sum_{i=1}^{N} |\mathcal{L}u_{\bXheta}(\bx^i_t, t) - f(\bx^i_t, t)|^2 \Big],
\end{equation}
where $\bx^i_t\sim\mu_t$.%, where the solution function is approximated with a parametrized neural network function $u_{\bXheta}(\bx, t)$ with all parameters being nonlinearly entered in the function. The integral above is taken with respect to some dynamical measure $\mu_t$, which will be discussed in detail in Section \ref{MSM}. 
%The selection of $\mu_t$ residual is defined as 
%\begin{equation}\label{residual}
%    r_t(\bx) = r(\bx, t) = \partial_t u(\bx, t) - \mathcal{L}u.
%\end{equation}

\begin{re}
    Of course, one can consider the time variable $t$ as an additional spatial variable and consider the PDE as just a stationary PDE. But in the scenario of efficiency by adaptive sampling, this may not be natural for (1) the PDE usually should be considered on some fixed time points (equal or non-equal time length); and (2) all time points should be computed with equal efficiency and the same sample size, otherwise some time points may even not be computed; and most importantly, (3) we may need to compute the PDE in infinite time and make prediction of the solution. 
\end{re}

Since we use Monte Carlo approximation of the loss \eqref{timePDENN}, we only need to consider the dynamics of the collocation samples. So we consider the flow of the adaptive sampling and set up the adaptive collocation point $\bx$ at time $t=0$ evolving according to $\bx_t=\bX_t(\bx)$ with the PDE solution, that is, a collocation point $\bx$ is evolved to a new point $\bX_t(\bx)$ and then used as a collocation point for learning. We now need to make $\bX_t(\bx)$ a flow that can fit our need for adaptivity. In addition, we don't need to sample again and again at different time $t$, we just need to map the initial sample $\bx$ at a time point $t$ with $\bX_t$. So $\bx$ can be viewed as the common parameterization for all time points $t\in[0,T]$ in moving mesh methods. 

The following theorem is well known for the transport equation, and will be restated to fit our need. Its proof is standard and given in Appendix~\ref{appendix: proof}. 

\begin{thm}\label{dynamics_thm2}
    Let $\bv_t:\mathbb{R}^n\rightarrow\mathbb{R}^n$ be a bounded and twice differentiable velocity field of the flow mapping $\bX_t$, that is, $\bX_t(\bx)$ follows the ODE system
    \begin{equation}\label{ODE_T_t}
        \frac{\partial \bX_t}{\partial t}(\bx) = \bv(t, \bX_t(\bx)).
    \end{equation}
    Suppose the density functions $p_t$ satisfies the transport equation
    \begin{equation}\label{PDE_v_t2}
        \frac{\partial p_t}{\partial t} + \nabla\cdot(p_t \bv_t) = 0.
    \end{equation}
    Then $p_t$ is the probability density of the probability flow mapping $\bX_t$, whose velocity field is $\bv_t$, in the sense that for all test functions $\varphi\in C^\infty_c (\mathbb{R}^n)$, 
    \begin{equation}\label{prob_p_t}
        \int_{\mathbb{R}^n}\varphi(\bx) p_t(\bx)\,d\bx = \int_{\mathbb{R}^n}\varphi(\bX_t(\bx)) p_0(\bx)\,d\bx.
    \end{equation}
    Conversely, if $p_t$ is the probability density evolution of $\bX_t$, whose velocity field is $\bv_t$, then it satisfies the transport equation \eqref{PDE_v_t2}.
\end{thm}

\subsection{The formulation of the moving sample method}\label{Subsec:MSM}
Now, the problem becomes how to determine the dynamics of $\bX_t$, namely, to find a velocity field $\bv_t$ such that numerical solution of \eqref{ODE_T_t} evolves the samples $\bX_t$ toward the distribution $\bX_t\sim p_t\propto r_t^2$. Then we obtain adaptive sampling flows of the collocation points and use them to solve \eqref{timePDENN} at any time $t$. %In \cite{Yuxiao_Coupling_2023}, this dynamics is just a coupled SDE that generates a Gibbs measure, where the SDE needs to be solved separately with a different step size, which is not efficient. 

Since $r_t^2$ is not normalized, we need to consider the transport equation for it. To simplify, we consider the transport equation for the logarithmic density $\log p_t$. Firstly, we have
\begin{align}
    \frac{\partial}{\partial t}\log p_t = & \frac{1}{p_t}\frac{\partial}{\partial t} p_t \notag\\
    = & - \frac{1}{p_t}\nabla\cdot(p_t\bv_t) \notag\\
    = & - \frac{1}{p_t}\nabla p_t\cdot\bv_t - \nabla\cdot\bv_t \notag\\
    = & - \nabla(\log p_t)\cdot\bv_t - \nabla\cdot\bv_t.\label{eqn_log_transport}
\end{align}
For the terms with $p_t$ in \eqref{eqn_log_transport}, we have
\begin{equation}\label{dt_log_density}
    \frac{\partial}{\partial t}\log p_t = \frac{\partial}{\partial t}\log\Big(\frac{r_t^2}{\int_{\mathbb{R}^n} r_t^2\,d\bx}\Big)
    = \frac{2}{r_t}\frac{\partial}{\partial t}r_t - \frac{1}{\int_{\mathbb{R}^n} r_t^2\,d\bx}\frac{d}{dt}\int_{\mathbb{R}^n} r_t^2\,d\bx,
\end{equation}
and
\begin{equation}\label{dx_log_density}
    \nabla\log p_t = \nabla\log\Big(\frac{r_t^2}{\int_{\mathbb{R}^n} r_t^2\,d\bx}\Big) = \frac{2}{r_t}\nabla r_t.
\end{equation}
For the first term in \eqref{dt_log_density} and the term in \eqref{dx_log_density}, we can just use the auto-differentiation to calculate exactly in every learning step. For the second term in \eqref{dt_log_density}, we will just use difference quotient
$$
\frac{d}{dt}\int_{\mathbb{R}^n} r_t^2\,d\bx \approx \frac{\int_{\mathbb{R}^n} r_{t_{i+1}}^2 - r_{t_i}^2\,d\bx}{t_{i+1}-t_i}
$$
to approximate it, in which the two terms in the numerator have been already calculated in the last iteration of the learning process. 

Combing the above equations, we have an equation for the velocity field. But this equation needs not to be solved exactly, since we only need to find an approximate sampling flow for solving the time-dependent PDE adaptively, like moving mesh method. So an adaptive method is not necessary to be applied for this part. In this work, we will also use PINNs to solve it.

The solution of either the transport equation \eqref{PDE_v_t2} or our formulation \eqref{eqn_log_transport} is not unique. We only need one of them to determine the dynamics of the sampling points efficiently. By the Helmholtz decomposition theorem \cite{arfken2011mathematical}, any differentiable vector field can be decomposed into a gradient field and a divergence-free vector field. We simply assume the vector field to be a pure gradient form $\bv_t = \nabla {\boldsymbol{\phi}}_t$ with a moving potential function ${\boldsymbol{\phi}}_t$. %Building on the established formulation of velocity field, we leverage the Helmholtz decomposition theorem \cite{arfken2011mathematical} to restrict the field to a pure gradient form $\bv_t = \nabla {\boldsymbol{\phi}}_t$. 
This restriction preserves the divergence structure relevant to the continuity constraint while eliminating non-identifiable solenoidal components that do not contribute to density evolution under the given boundary conditions. By reducing the parameter space and enforcing an irrotational flow, the gradient formulation improves numerical stability, and aligns the learned dynamics with the underlying physical and probabilistic structure of the problem.

\section{Algorithm Implementation}\label{Algorithm Implementation}
In the previous section, we presented the equations satisfied by $p_t$ and velocity field $\bv_t$, along with the differential equation used to solve for $\bv_t$. In this section, we describe the implementation of the algorithm.

We employ two fully connected neural networks $u_{\bXheta}(\bx, t), {\boldsymbol{\phi}}_{\boldsymbol{\eta}}(\bx, t)$ to separately approximate the solution $u$ of the partial differential equation and the potential function of the grid velocity field $\bv_t$, where ${\boldsymbol{\phi}}_{\boldsymbol{\eta}}$ is a scalar potential network, $\bXheta, {\boldsymbol{\eta}}$ denote the parameters of the neural networks. Define $\bv_{\boldsymbol{\eta}} = \nabla {\boldsymbol{\phi}}_{\boldsymbol{\eta}}$, then $\nabla \times \bv_{\boldsymbol{\eta}} = \nabla \times(\nabla {\boldsymbol{\phi}}_{\boldsymbol{\eta}}) = 0$, ensuring that $\bv_{\boldsymbol{\eta}}$ is an irrotational velocity field. From %\cite{raissi2019physics} and 
\eqref{eqn_log_transport}, \eqref{dt_log_density}, and \eqref{dx_log_density}, the loss function of $u_{\bXheta}, \bv_{\boldsymbol{\eta}}$ is then defined separately as 
\begin{align}
    \text{Loss}_u(\bXheta) = & \int_0^T\int_\Omega |\mathcal{L}u_{\bXheta} - f|^2 \,\mu_t(d\bx)dt + \beta_1 \int_{\Omega} |\mathbf{\mathcal{I}}u_{\bXheta} - \bu_0|^2 \mu_0(d\bx) \notag\\
    & + \beta_2 \int_0^T\int_{\partial \Omega} |\mathcal{B}u_{\bXheta} - u_b|^2 \,\mu_b(d\bx)dt, \label{loss_u}
\end{align}
\begin{align}
    \text{Loss}_{\bv}({\boldsymbol{\eta}}) &= \int_0^T\int_\Omega (2\frac{\partial}{\partial t}r_t + 2\nabla r_t \cdot \bv_{\boldsymbol{\eta}} + r_t \nabla\cdot \bv_{\boldsymbol{\eta}} - \frac{r_t}{\int_\Omega r_t^2 d\bx} R_t)^2 \,d\bx dt \notag\\
    &= \int_0^T\int_\Omega (2\frac{\partial}{\partial t}r_t + 2\nabla r_t \cdot \nabla {\boldsymbol{\phi}}_{\boldsymbol{\eta}} + r_t \Delta {\boldsymbol{\phi}}_{\boldsymbol{\eta}} - \frac{r_t}{\int_\Omega r_t^2 d\bx} R_t)^2 \,d\bx dt, \label{loss_v}
\end{align}
where $\beta_1, \beta_2$ are used to balance the components of the loss function, and
$$
R_t=\frac{d}{d t} \int_\Omega r_t^2 \,d\bx.
$$

A sketch of the overall algorithm is presented in Alg. \ref{alg:MSM}.
\begin{algorithm}
\caption{Moving sample method\label{alg:MSM}}
\begin{algorithmic}[1]
    \REQUIRE Initial networks $u_{\bXheta}, \boldsymbol{\phi}_{\boldsymbol{\eta}}$, the collocation points for PDE loss $S = \{(\bx^{(i)}_{\Omega}, t^{(i)})\}_{i=1}^{N\times N_t}$, training set for initial condition $S_0 = \{\bx_0^{(i)}\}_{i=1}^{N_0}$,  training set for boundary condition $S_{bdry} = \{(\bx_{\partial \Omega}^{(i)}, t^{(i)})\}_{i=1}^{N_b\times N_t}$, the number of training iteration $M$, training epochs $M_1, M_2$ for $u_{\bXheta}, \boldsymbol{\phi}_{\boldsymbol{\eta}}$ in each iteration.
    \FOR{$i = 1, 2, \dots, M$}
        \STATE Train the model $u_{\bXheta}$ for $M_1$ epochs by descending the stochastic gradient of $\text{Loss}_u$ \eqref{loss_u} on $S, S_0, S_{bdry}$.
        \STATE Train the model $\boldsymbol{\phi}_{\boldsymbol{\eta}}$ for $M_2$ epochs by descending the stochastic gradient of $\text{Loss}_{\bv}$ \eqref{loss_v} on $S$.
        \STATE Sample new $\bX_0$ according to the initial condition $\mathcal{I}u(\bx, 0)$.
        \STATE Compute new $\bX_t$ from $\bv_{\boldsymbol{\eta}}$ and $\bX_0$ numerically from \eqref{ODE_T_t} to get the new adaptive samples $S_{new} = \{(\bx_{adaptive}^{(i)}, t^{(i)})\}_{i=1}^{N_1\times N_t}$.
        \STATE Update the training set $S$ by $S = S\cup S_{new}$.
    \ENDFOR
    \STATE Train the model $u_{\bXheta}$ for $M_{final}$ epochs by descending the stochastic gradient of $\text{Loss}_u$ \eqref{loss_u} on $S, S_0, S_{bdry}$.
    \ENSURE $u_{\bXheta}, S$.
\end{algorithmic}
\end{algorithm}

Here we would like to specifically point out that we do not use any uniqueness condition for the training of the velocity field $\bv_{\boldsymbol{\eta}} = \nabla {\boldsymbol{\phi}}_{\boldsymbol{\eta}}$, since we just need one velocity field that drives the samples to move according to the residual changes. However, this may course the samples to move out of the computational region. For example, in Experiment \ref{burgers}, the samples are on the interface of the wave, but some will move out of the region. But this does not affect the training of the examples in our work and most practical cases, since the boundary conditions, as well as the corresponding region, are just conditions to determine the uniqueness of the solution, but the solution is not really restricted in this region. In other words, the samples out of the region can be still used as learning collocations. When it is necessary to restrict the solution and the training samples in a specific region, we can add conditions to the potential ${\boldsymbol{\phi}}_{\boldsymbol{\eta}}$ of the velocity field. 

In the sampling component, we choose $\mu_0$ to be related to the initial condition of the equation. In the numerical experiments, we employed two different sampling strategies to sample initial training set $S_0 = \{\bx_0^{(i)}\}_{i=1}^{N_0}$. The first strategy distributes the sampling points in proportion to $u_0$, while the second distributes them in proportion to the squared Euclidean norm of the gradient of $u_0$. Both strategies also include adding a number of uniformly distributed sampling points. For boundary training set $S_{bdry}$ and PDE training set $S$, we uniformly divide the time interval $[0, T]$ into $N_t$ segments, and compute the PDE loss and boundary loss at these discrete time points. For the boundary conditions, we take $\mu_b$ as a uniform measure, randomly generate $N_b$ samples on $\partial \Omega$, and compute the boundary loss at these points across all $N_t$ time steps. The boundary training set is denoted by $S_{bdry}=\{(\bx_{\partial \Omega}^{(i)}, t^{(i)})\}_{i=1}^{N_b\times N_t}$. Similarly, the training set $S$ for computing the PDE loss consists of a small subset of uniformly distributed sampling points and adaptive sampling points. The uniformly distributed points remain identical across all $N_t$ time steps, whereas the adaptive sampling points move under the influence of the velocity field.

\begin{re}
    As mentioned at the beginning of Subsection \ref{dynamics}, the sampling distribution $p_t$ used in Theorem \ref{dynamics_thm2} can be chosen to be the normalization of $r_t^{2\gamma}$, $|u|^{\gamma}$, $|\nabla_x u|^{\gamma}$, or their mixture, as in the moving mesh methods. But no matter which we choose, it is an explicit formulation of the density function of this distribution by using the approximation $u_{\bXheta}$ obtained in the last iteration of PINNs learning. Furthermore, no matter how the sampling distribution in the first iteration is chosen (e.g. uniform distribution in vanilla PINNs), and thus how good the first approximation is, the above explicit formulated density can always describe the current demand of adaptive samples. 

    However, there is still one computational perspective we should consider. For the same reason of introducing adaptive sampling, if $p_t$ is highly concentrated due to high singularity of the PDE, it is still not easy to train the velocity field $\bv_t$. To be specific, uniformly distributed samples may not capture high concentration of $p_t$, and thus may give a trivial velocity field. So, even if it is not necessary to get $\bv_t$ exactly, we still need an effective approximation to guide the movement of adaptive samples. In the first iteration for training $\bv_t$, we observe that samples according to the singularity of the initial condition is already enough for training $\bv_t$ well. This means in the training of $\bv_t$, we can include $S_0$, at least the samples for $t_0$. After the first iteration, we obtain samples that can capture concentration of $p_t$. Then we can use the samples from the previous iterations. 
\end{re}

\section{Numerical Experiments}\label{Numerical Experiments}
In this section, we conduct a series of numerical experiments to test the effectiveness of this method. To quantitatively evaluate the accuracy of the results, we use the relative $L_2$ error $\frac{\|u_{\bXheta} - u^*\|_2}{\|u^*\|_2}$ to measure the global approximation performance, and the $L_\infty$ error $\max|u_{\bXheta} - u^*|$ to indicate the ability of capturing the singularity of the approximation, for Experiments~\ref{Allen-Cahn}--\ref{fokker planck}, where $u^*(\bx, t)$ is the exact solution. For Experiment~\ref{high-dim advection}, we employ a weighted $L_2$ norm and $L_\infty$ norm, with the weights concentrated in high-density regions. In all experiments, we set the networks for both $u_{\bXheta}$ and $\boldsymbol{\phi}_{\boldsymbol{\eta}}$ to have $3$ hidden layers, with $64$ neurons per hidden layer. The activation functions are set to {\it tanh} for all hidden layers, while no activation function is applied to the output layer. The neural networks are trained using the Adam optimizer with a fixed learning rate of $0.001$, and the training is carried out for a total of $M=5$ iterations. The number of training epochs for $u_{\bXheta}$ varies across experiments: in Experiments~\ref{Allen-Cahn} and~\ref{fokker planck}, each iteration uses $M_1 = 6000$ epochs, while in the remaining experiments $1500$ epochs are used. The velocity potential field network $\boldsymbol{\phi}_{\boldsymbol{\eta}}$ is trained with $M_2 = 1000$ epochs in all cases. A tailored number of adaptive sampling points is employed for the initial condition and PDE components, together with a fixed set of uniformly distributed points. The specific sampling numbers for each example are provided in the corresponding subsections.

\begin{table}[htbp]
\centering
\begin{tabular}{lcccc}
\hline
 & \textbf{Method} & \textbf{Relative $L_2$ Error} & \textbf{$L_\infty$ Error} \\
\hline
Allen-Cahn & PINNs & $9.829\times10^{-3}$ & $2.011\times10^{-1}$ \\
            & MSM-PINNs & $3.597\times10^{-3}$ & $3.903\times10^{-2}$ \\
\hline
Rotation & PINNs & $8.967\times10^{-2}$ & $1.370\times10^{-1}$ \\
            & MSM-PINNs & $5.993\times10^{-2}$ & $3.020\times10^{-2}$ \\
\hline
Burgers' & PINNs & $6.539\times10^{-3}$ & $1.734\times10^{-1}$ \\
            & MSM-PINNs & $1.784\times10^{-3}$ & $5.387\times10^{-2}$ \\
\hline
Fokker-Planck & PINNs & $3.637\times10^{-2}$ & $2.498\times10^{-1}$ \\
            & MSM-PINNs & $1.761\times10^{-2}$ & $9.167\times10^{-2}$ \\
\hline
High-Dimensional  & PINNs & - & - \\
           Advection & MSM-PINNs & $1.177\times10^{-2}$ & $1.241\times10^{-1}$ \\
\hline
\end{tabular}
\caption{Comparison of relative $L_2$ and $L_\infty$ errors for different models with vanilla PINNs and PINNs with MSM. In each experiment, the results are obtained with the same number of sampling points, training epochs, and network architecture. In particular, for the advection case \ref{high-dim advection}, the values are weighted errors. The absence of error indicates that the method fails to produce a meaningful approximation.}
\label{tab:error}
\end{table}

\subsection{Allen--Cahn equation}\label{Allen-Cahn}
We first consider the following one-dimensional Allen--Cahn equation:
\begin{equation}\label{eq:AC}
\left\{
\begin{array}{ll}
    \frac{\partial u}{\partial t} = \alpha \frac{\partial^2 u}{\partial x^2} - \beta(u^3 - u), & x\in [-1, 1], t \in [0, 1], \\
    u(x, 0) = x^2 cos(\pi x), & x\in [-1, 1],\\
    u(-1, t) = u(1, t) = -1, & t \in [0, 1].
\end{array}\right.
\end{equation}
We set the coefficient $\alpha$ to $0.001$, $\beta$ to $5$ in this case.

This equation describes phase separation dynamics and the evolution of interfaces in a binary medium, which exhibits sharp transition layers near $x=0$ and $x=\pm 0.5$. We impose a hard-constraint transformation on the network output $y$ by defining $u_{\bXheta} = t(x^2 - 1)y + u_0$, so that the result inherently satisfies both the boundary and initial conditions, thereby improving training performance. We choose a uniform distribution as the sampling distribution for $\bX_0$ because the boundary and initial conditions in this example exhibit low singularity. Moreover, to further enhance the training performance, we sample the initial collocations $S$ uniformly distributed over the entire domain $\Omega \times [0, T]$, i.e., sampled jointly in time rather than solely at the $N_t$ discrete temporal grid points.

\begin{figure}
    \centering
    \includegraphics[width=0.79\linewidth]{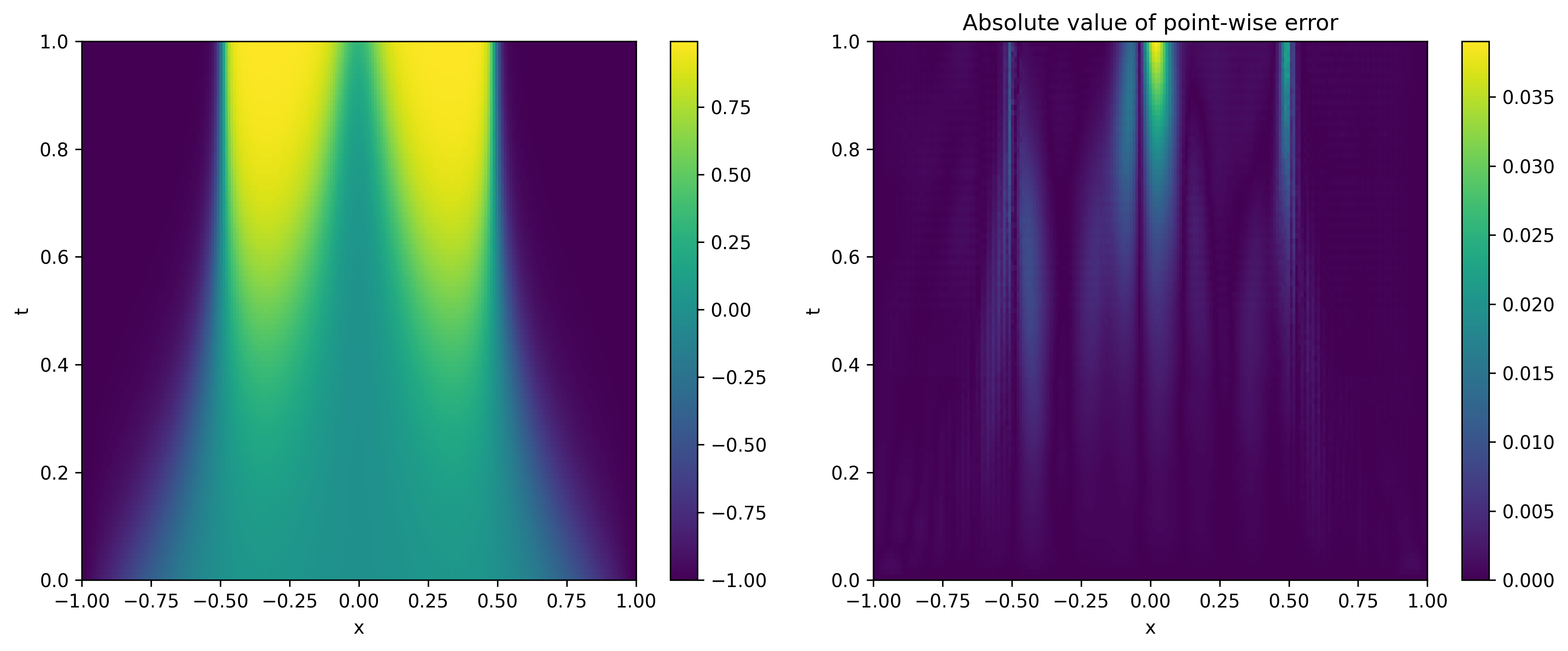}
    \caption{MSM-PINNs result (left) and the corresponding absolute error (right) for the Allen-Cahn equation in Experiment~\ref{Allen-Cahn}.}
    \label{fig:Allen-Cahn solution}
\end{figure}

\begin{figure}
    \centering
    \includegraphics[width=0.8\linewidth]{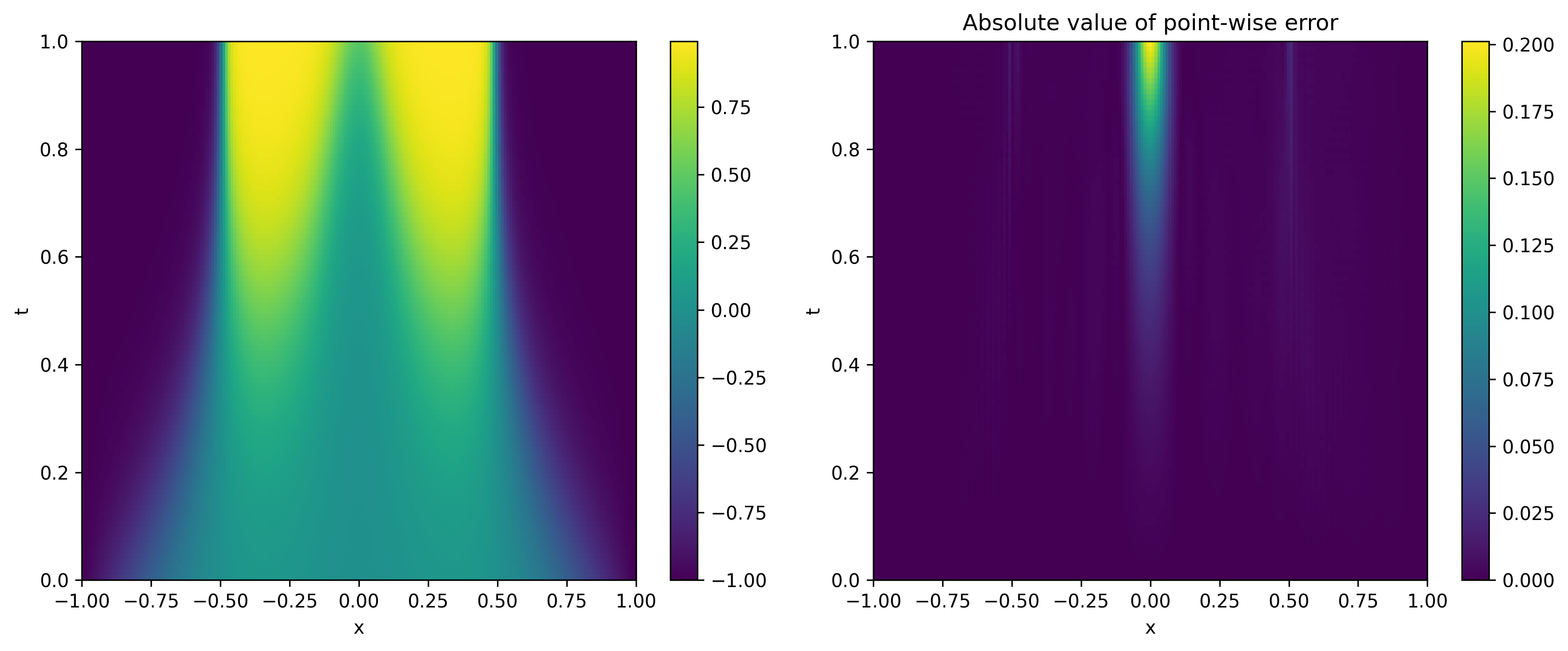}
    \caption{PINNs result (left) and the corresponding absolute error (right) for the Allen-Cahn equation in Experiment~\ref{Allen-Cahn}.}
    \label{fig:Allen-Cahn PINNs}
\end{figure}

\begin{figure}
    \centering
    \includegraphics[width=0.95\linewidth]{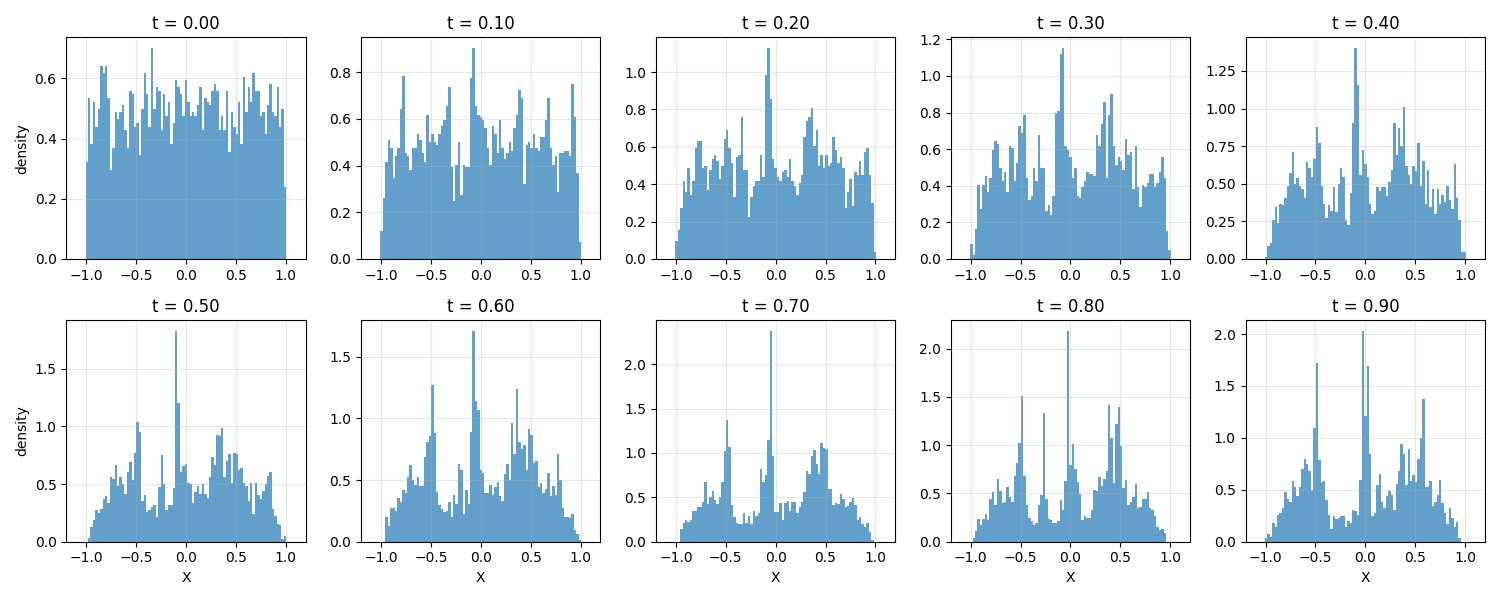}
    \caption{The approximate distribution of adaptive sampling points for the Allen-Cahn equation in Experiment~\ref{Allen-Cahn}.}
    \label{fig:Allen-Cahn sample movement}
\end{figure}

In this example, we set the number of uniform points $N=200$ and the number of adaptive moving samples added in each iteration $N_1=600$ in MSM-PINNs. For the training of PINNs, we use a total of $3200$ uniformly sampled collocations. It is the same to the total number of samples we used in MSM-PINNs. The numerical solutions and their corresponding absolute errors are shown in Fig.~\ref{fig:Allen-Cahn solution} and Fig.~\ref{fig:Allen-Cahn PINNs}. The evolution of the sampling point distribution generated by the proposed method is illustrated in Fig.~\ref{fig:Allen-Cahn sample movement}. Starting from a uniform distribution at $t=0$, the sampling points progressively migrate toward regions with stronger singularities, particularly around $x=0$ and $x=\pm 0.5$. This movement reflects the residual-driven velocity field learned by the model and demonstrates that the adaptive strategy effectively concentrates computational effort in high-residual regions.

A quantitative comparison of the numerical accuracy is provided in Table~\ref{tab:error}. For the Allen--Cahn equation, the proposed MSM-PINNs method achieves approximately $63\%$ improvement in the relative $L_2$ error and $81\%$ improvement (from $2.011 \times 10^{-1}$ to $3.903 \times 10^{-2}$, also see Fig.~\ref{fig:Allen-Cahn solution} and \ref{fig:Allen-Cahn PINNs}) in the $L_\infty$ error over standard PINNs. These results confirm that concentrating sampling points in high-residual regions significantly enhances the solution accuracy. In particular, the relative $L_2$ error only indicates the global performance improvement of the MSM-PINNs, whereas the $L_\infty$ error provides a clearer indication of the local behavior of the solution. Specifically, MSM-PINNs can capture the singularity around $x=0$ where PINNs could not approximate well (see the right figure in Fig.~\ref{fig:Allen-Cahn PINNs}), and then allocates more learning collocations there (see Fig.~\ref{fig:Allen-Cahn sample movement}) to improve the approximation of the solution singularity (see Fig.~\ref{fig:Allen-Cahn solution}). 

\subsection{Rotation equation}\label{rotation}

We next consider the following two-dimensional equation, which represents the rotation of a Gaussian kernel moving on the unit circle:
\begin{equation}\label{eq:movingGauss}
\left\{
\begin{array}{ll}
    \frac{\partial u}{\partial t} = \frac{\partial u}{\partial x}\sin t - \frac{\partial u}{\partial y}\cos t, & \text{in}\ \Omega\times[0, 1], \\
    u(x, y, 0) = e^{-\frac{1}{\alpha}((x - 1)^2 + y^2)}, & \text{in}\ \Omega,\\
    u(x, y, t) = e^{-\frac{1}{\alpha}((x - \cos t)^2 + (y - \sin t)^2)}, & \text{on}\ \partial\Omega \times [0, 1].
\end{array}\right.
\end{equation}
The exact solution is
$$    
u(x, y, t) = e^{-\frac{1}{\alpha}((x - \cos t)^2 + (y - \sin t)^2)}.
$$
In this case, we set $\alpha = 0.01$, $\Omega = (-0.2,1.2)\times(-0.2,1.2)$.

\begin{figure}
    \centering
    \includegraphics[width=0.95\linewidth]{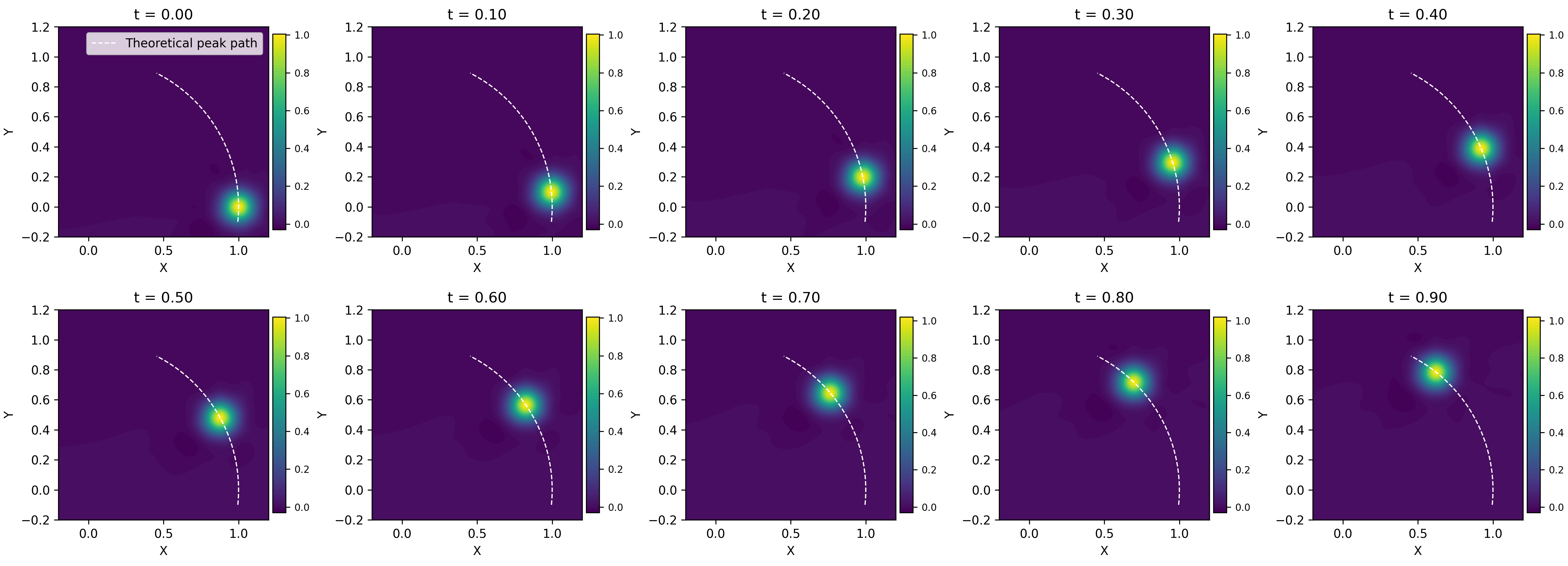}
    \caption{MSM-PINNs result for the rotation equation in Experiment~\ref{rotation}.}
    \label{fig:rotation solution}
\end{figure}

\begin{figure}
    \centering
    \includegraphics[width=0.95\linewidth]{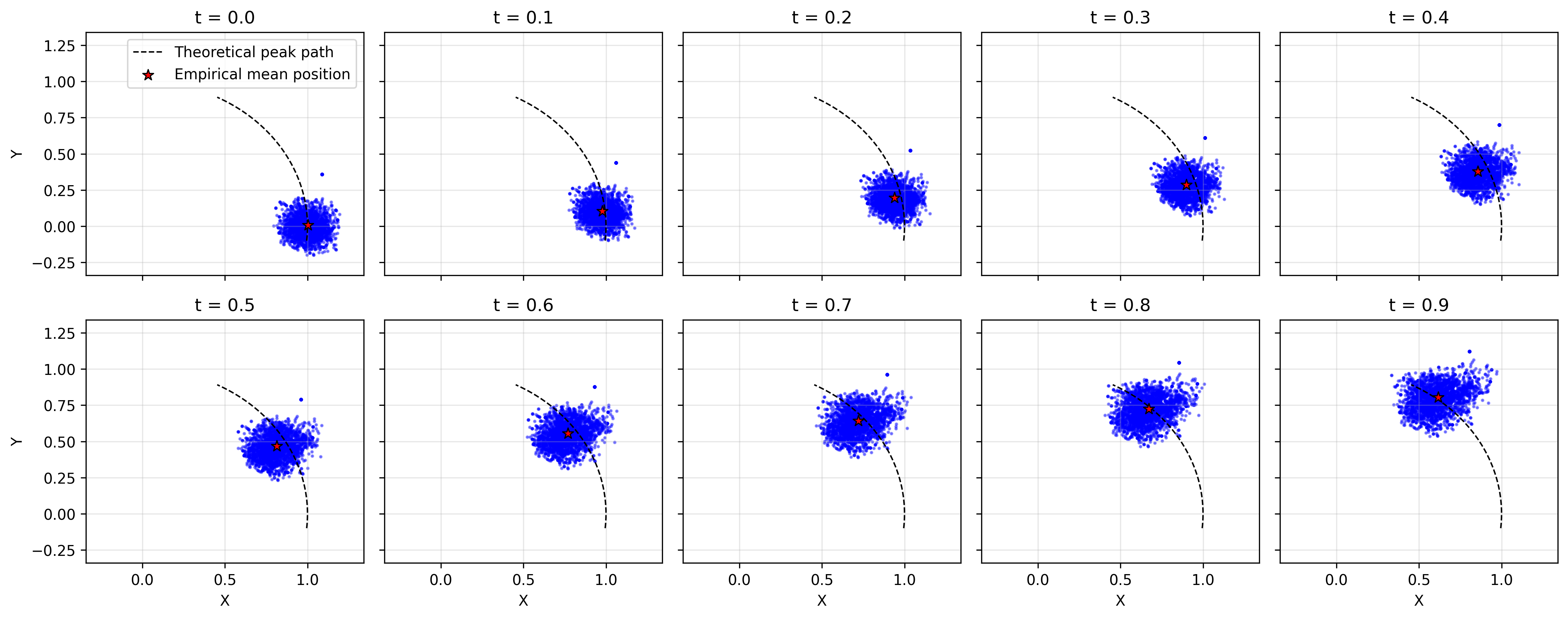}
    \caption{The trajectories of the adaptive sampling points for the rotation equation in Experiment~\ref{rotation}.}
    \label{fig:rotation sample movement}
\end{figure}

\begin{figure}
    \centering
    \includegraphics[width=0.95\linewidth]{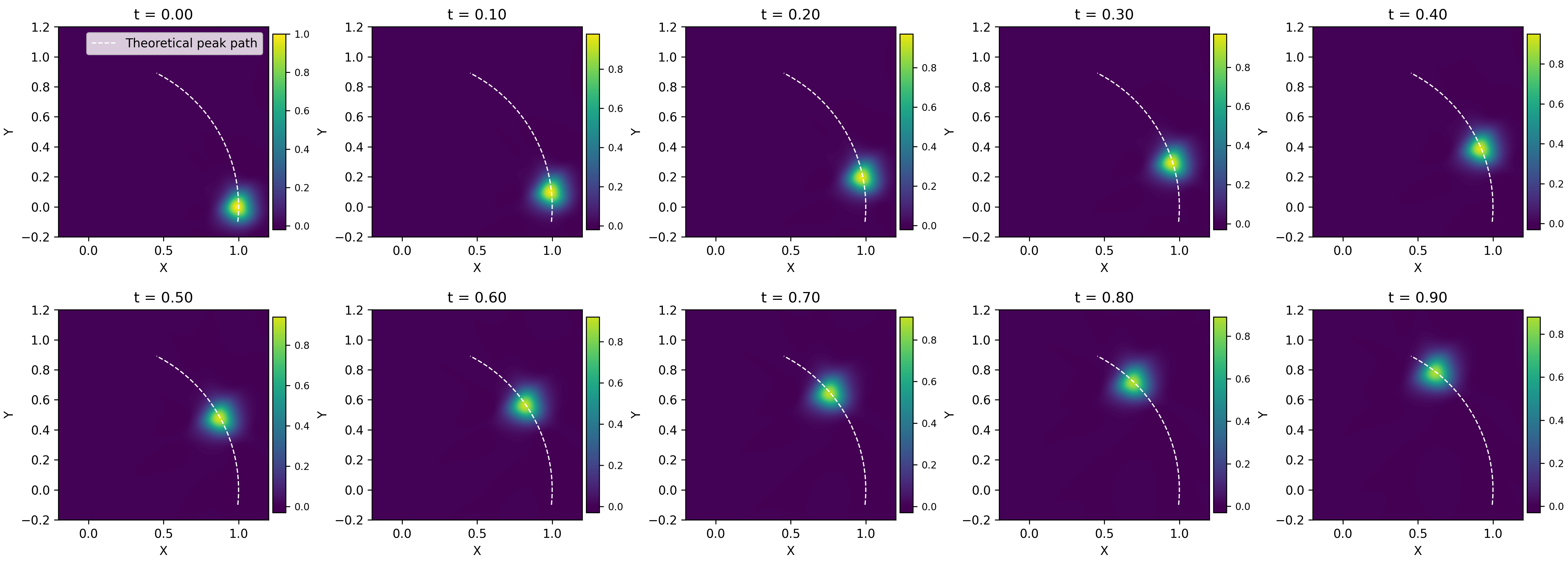}
    \caption{PINNs solution for the rotation equation in Experiment~\ref{rotation}.}
    \label{fig:rotation PINNs solution}
\end{figure}

\begin{figure}[htbp]
    \centering
    \includegraphics[width=0.95\linewidth]{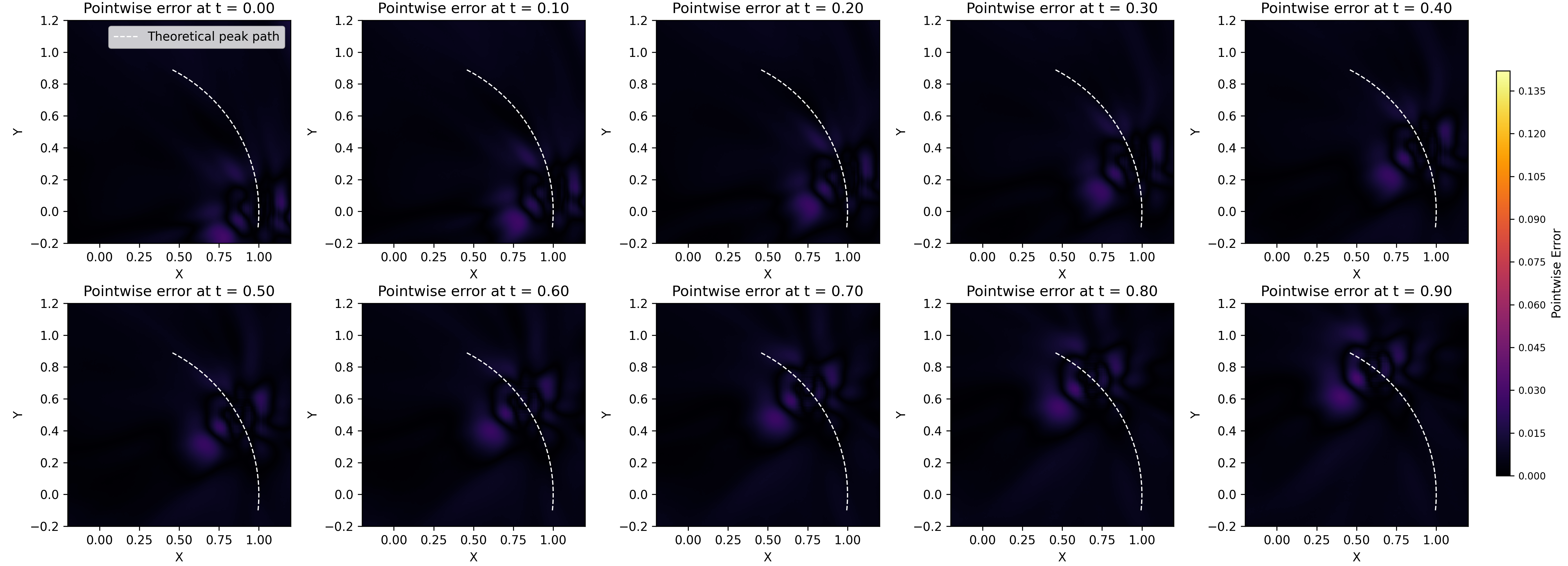}
    \par\vspace{0.5cm}    
    \includegraphics[width=0.95\linewidth]{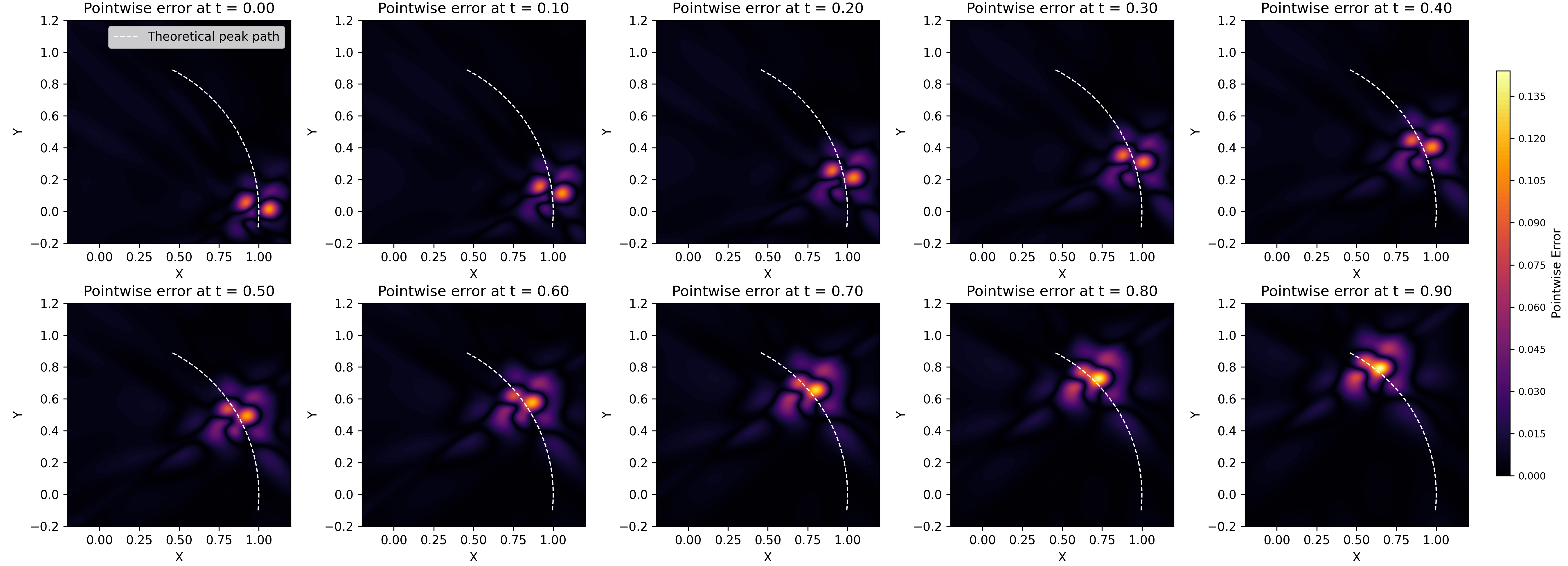}
    \caption{Point-wise errors comparison between results with MSM-PINNs (top) and PINNs (bottom).}
    \label{fig:rotation pointwise error}
\end{figure}

This equation models the passive transport of a distribution with high singularity under a prescribed rotational flow to show non-trivial (e.g., nonlinear) movement of adaptive samples. For this equation, we sample the initial training set $S_0$ proportional to the initial condition $u_0$ (one can also use other $|\nabla u_0|^\gamma$, just as the moving sample distributions). We deliberately set a smaller number of sample points to demonstrate the capability of our method to train the entire model with a smaller amount of data. More specifically, we set $N=1000, N_0=500, N_b=400$, and introduce $N_1 = 1000$ adaptive points after the first iteration, followed by $N_1 = 300$ points in each subsequent iteration. The numerical results are shown in Fig.~\ref{fig:rotation solution}. The trajectories of all adaptive sampling points, visualized in Fig.~\ref{fig:rotation sample movement}, reveal that most points follow the rotating wave peak closely. This behavior indicates that the residual-driven velocity field successfully captures the underlying rotational dynamics and its singularity, enabling the sampling points to remain aligned with the moving high-gradient region throughout the evolution.

Across different training iterations, the spatial distribution of adaptive sampling points exhibits characteristic variations. While most points consistently track the rotating wave crest, a subset gradually diverges toward other regions of the domain (see Fig.~\ref{fig:rotation sample movement}). This divergence typically appears in the later iterations, when the residual around the main peak has been substantially reduced, prompting the sampling strategy to shift attention to secondary regions where the residual remains non-negligible. In this example, the effect manifests as a mild dispersion of sampling trajectories, whereas in Example~\ref{high-dim advection} it leads to slight changes in movement direction or the velocity of points generated at different iterations. Such behavior reflects the flexibility of the adaptive method and contributes to the improved accuracy of the proposed method.

For the training of PINNs, we use a total of $2200$ uniformly sampled collocations for the PDE loss, which is the same as the total number of samples used in MSM-PINNs. A representative result obtained using standard PINNs is shown in Fig.~\ref{fig:rotation PINNs solution}, and Fig.~\ref{fig:rotation pointwise error} presents a comparison of the point-wise errors between MSM-PINNs and PINNs under the same scaling. Compared with the proposed method, the PINNs solution exhibits noticeable instability and reduced accuracy near the rotating peak, as summarized in Table~\ref{tab:error}. It can be seen that PINNs fail to accurately capture the highly singular peak, exhibiting instability near the peak at all time instances, particularly when $t$ becomes large (see Fig.~\ref{fig:rotation PINNs solution} and Fig.~\ref{fig:rotation pointwise error}). For larger values of $t$, the PINNs solution often suffers from dissipation and gradually loses the structural integrity of the rotating Gaussian. Further details of such results are provided in the Appendix~\ref{appendix: PINNs}, illustrating the instability of the PINNs approach. In contrast, the proposed method remains stable even with relatively few sampling points, as the adaptive strategy continuously reallocates points to the moving high-residual region and prevents dissipation.

\subsection{Burgers’ equation}\label{burgers}

Consider the following Burgers' equation:
\begin{equation}\label{eq:Burgers}
\left\{
\begin{array}{ll}
    \frac{\partial u}{\partial t} = \alpha (\frac{\partial^2 u}{\partial x^2} + \frac{\partial^2 u}{\partial y^2}) - u(\frac{\partial u}{\partial x} + \frac{\partial u}{\partial y}), & \text{in}\ \Omega\times[0, 1], \\
    u(x, y, 0) = \frac{1}{1+e^{\frac{x+y}{2\alpha}}},  & \text{in}\ \Omega,\\
    u(x, y, t) = \frac{1}{1+e^{\frac{x+y-t}{2\alpha}}}, & \text{on}\ \partial\Omega\times[0, 1].
\end{array}\right.
\end{equation}
The exact solution is 
$$
u(x, y, t) = \frac{1}{1+e^{\frac{x+y-t}{2\alpha}}}.
$$
In this case, we set $\alpha=0.001$, $\Omega=(-1,1)\times(-1,1)$.

\begin{figure}
    \centering
    \includegraphics[width=0.95\linewidth]{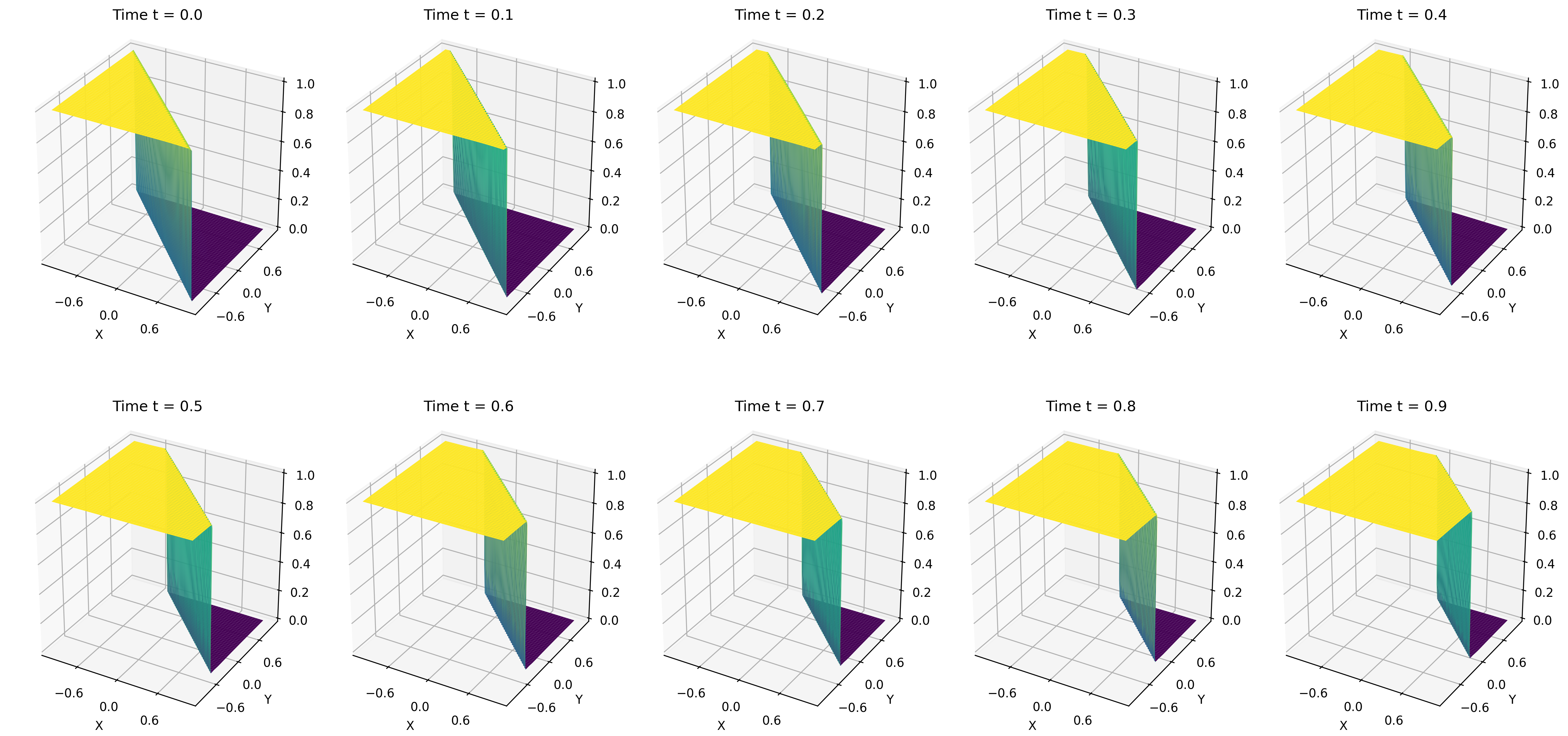}
    \caption{MSM-PINNs result for the Burgers' equation in Experiment~\ref{burgers}.}
    \label{fig:burgers solution}
\end{figure}

\begin{figure}
    \centering
    \includegraphics[width=0.95\linewidth]{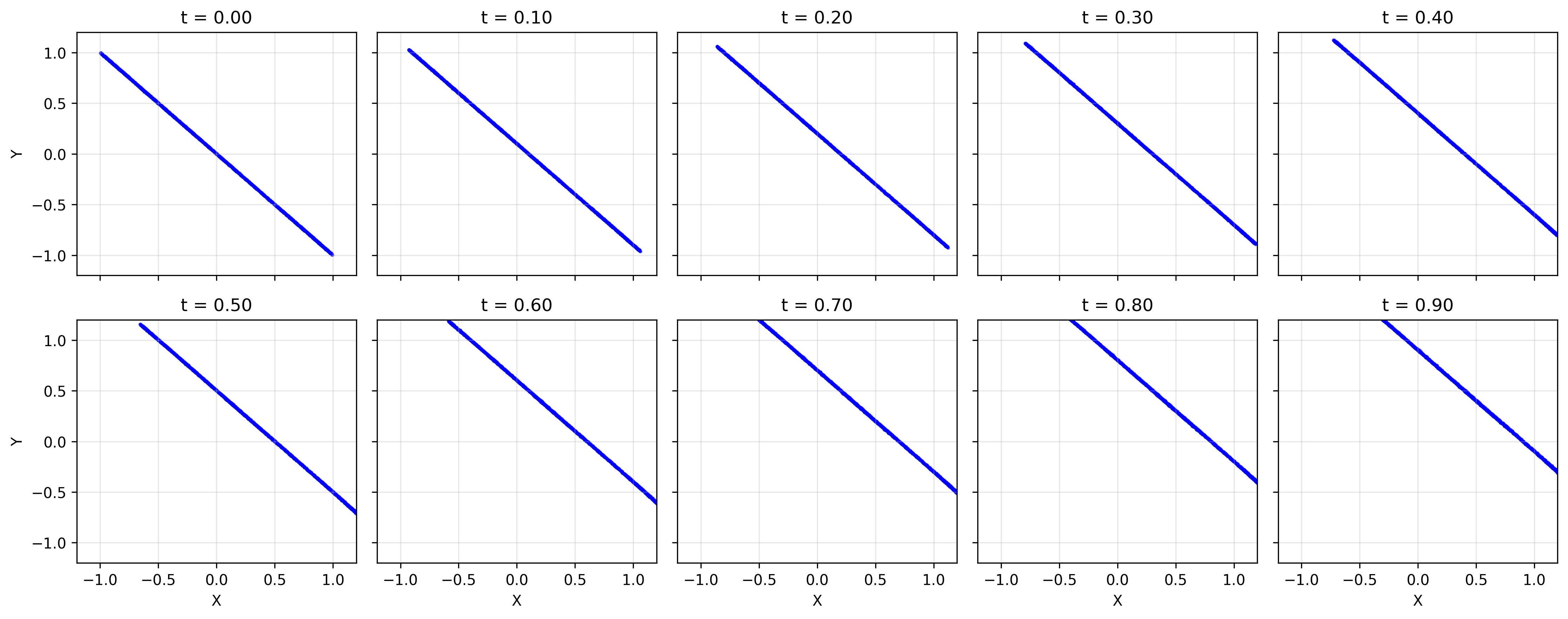}
    \caption{The trajectories of the adaptive sampling points for the Burgers' equation in Experiment~\ref{burgers}.}
    \label{fig:burgers sample movement}
\end{figure}

\begin{figure}
    \centering
    \includegraphics[width=0.95\linewidth]{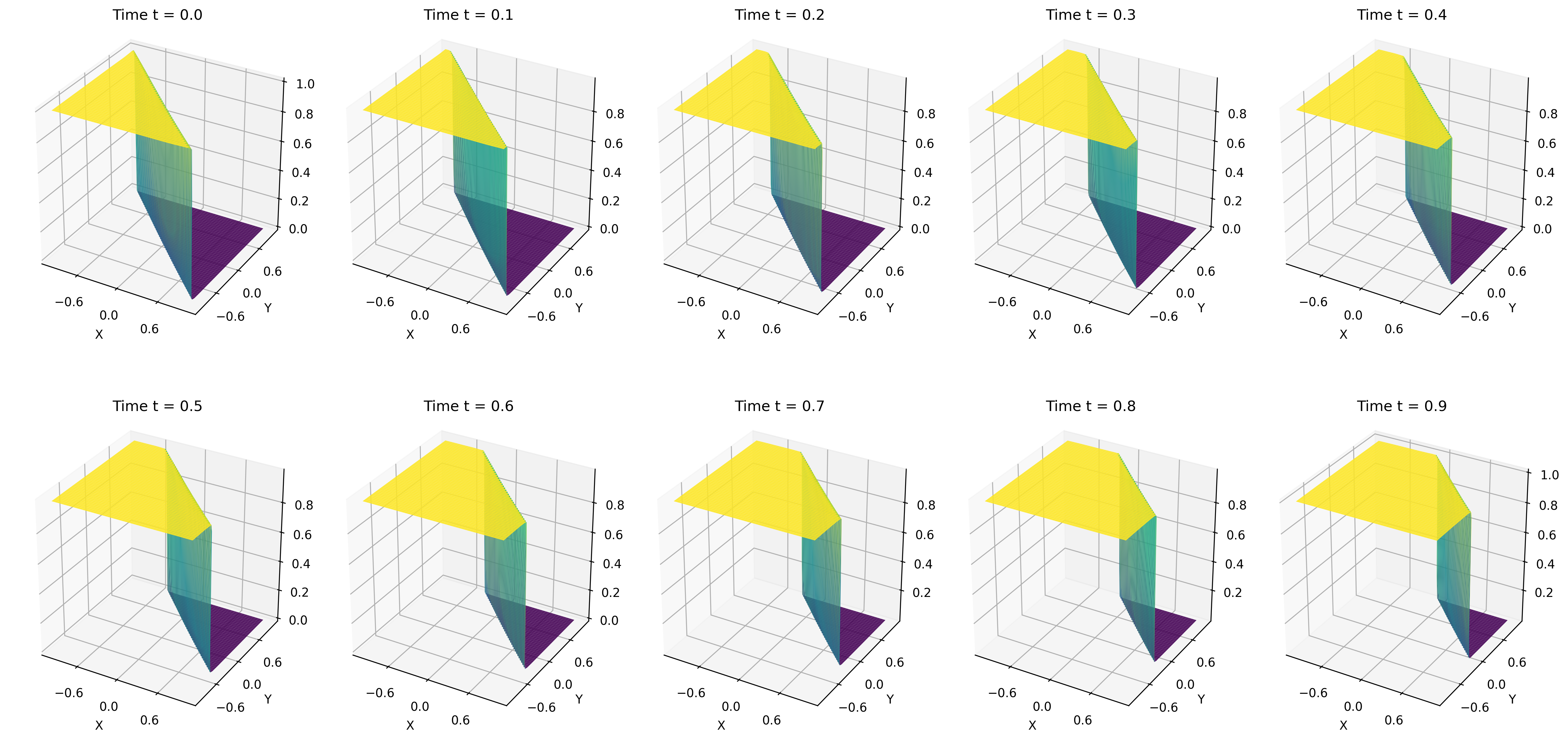}
    \caption{PINNs solution for the Burgers' equation in Experiment~\ref{burgers}.}
    \label{fig:burgers PINNs}
\end{figure}

This viscous Burgers' equation describes the evolution of a scalar field under diffusion with viscosity $\alpha$, which exhibits pronounced singularity at the interface when $\alpha$ is small. For this equation, we sample the initial training set $S_0$ proportional to the squared gradient of initial condition $u_0$. To further test the robustness of the proposed method, we deliberately use a relatively small number of sampling points. Specifically, we set $N=1200, N_0=500, N_b=200$, and introduce $N_1 = 1200$ adaptive points after the first iteration, followed by $N_1 = 300$ points in each subsequent iteration. The solution obtained by our method, together with the trajectories of the adaptive sampling points, is shown in Fig.~\ref{fig:burgers solution} and Fig.~\ref{fig:burgers sample movement}. The sampling points move together with the interface and remain concentrated around it as it propagates, which helps the network capture the sharp transition accurately.

\begin{figure}
\centering
\begin{subfigure}{0.32\textwidth}
    \centering
    \includegraphics[width=\linewidth]{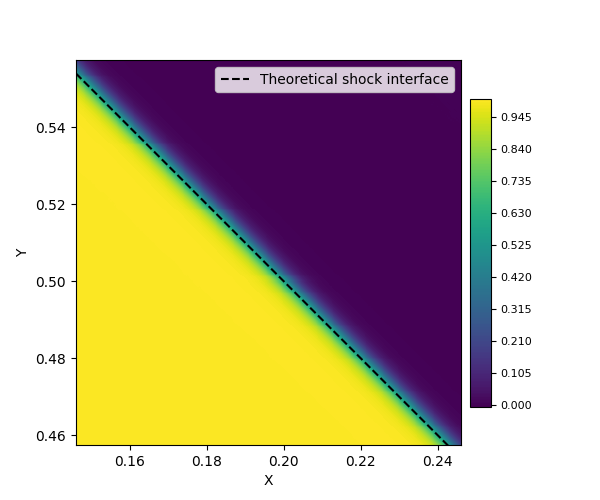}
\end{subfigure}
\begin{subfigure}{0.32\textwidth}
    \centering
    \includegraphics[width=\linewidth]{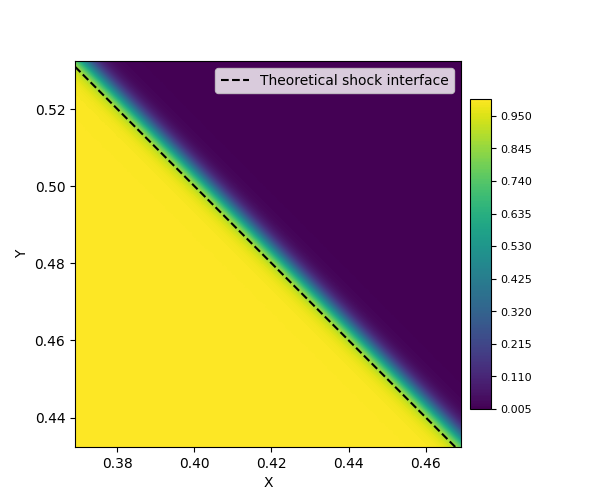}
\end{subfigure}
\begin{subfigure}{0.32\textwidth}
    \centering
    \includegraphics[width=\linewidth]{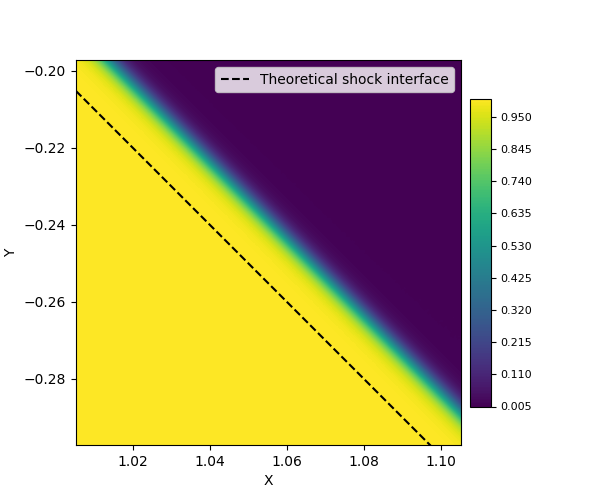}
\end{subfigure}
\caption{Enlarged views around the regions of maximum absolute error for the Burgers' equation in Experiment~\ref{burgers}: (left) MSM-PINNs; (middle and right) PINNs.}
\label{fig:burgers_local}
\end{figure}

For the training of PINNs, we use a total of $2400$ uniformly sampled collocations for the PDE loss, which is the same as the total number of samples used in MSM-PINNs. The results produced by standard PINNs under the same settings are unstable. As illustrated in Fig.~\ref{fig:burgers PINNs}, the solution shown represents the best result obtained among dozens of PINNs runs, which attains a relative $L_2$ error of $6.539\times10^{-3}$. However, this solution exhibits noticeable errors in singular regions. The enlarged views in Fig.~\ref{fig:burgers_local} highlight this issue: although the overall shape of the solution appears correct, the wave peak predicted by PINNs is noticeably shifted from its true position. The detailed error comparison is provided in Table \ref{tab:error}, where the reported $L_\infty$ error corresponds to the best outcome across multiple runs. Nevertheless, this value is not essentially different from other runs with $L_\infty$ error close to $1$ (the $L_\infty$ norm of the exact solution is $1$), all of which consistently reflect the difficulty PINNs face in accurately determining the interface location. The right figure in Fig.~\ref{fig:burgers_local} shows a usual example result seemingly correct, but the interface is actually away from that of the exact solution. In many runs, the deviation is even more pronounced, as shown by additional examples in Appendix~\ref{appendix: PINNs}. The middle one in Fig.~\ref{fig:burgers_local} shows the best result selected from dozens of selections, yet still not as correct as the left figure obtained with just a single random run of MSM-PINNs. It is important to notice that we can select the best one from the PINNs results since we know the exact solution. 

In contrast, the proposed method maintains stable performance across iterations. By continuously moving sampling points to the singular interface, it avoids the peak-shift phenomenon observed in PINNs and achieves higher accuracy with fewer sampling points. And any single random run of MSM-PINNs can get a similar approximation as shown in Figure \ref{fig:burgers solution} and similar error shown in Table \ref{tab:error}.

\subsection{Fokker-Planck equation}\label{fokker planck}

The Fokker-Planck equation (FPE) \cite{risken1989fokker} is a partial differential equation that describes how the probability density function of a stochastic process evolves over time under the influence of drift and diffusion. Consider the following stochastic differential equation:
\begin{equation}\label{eq:SDE}
    \left\{
    \begin{aligned}
        dx &= (-4(x-e^{-t})((x-e^{-t})^2 + (y-e^{-t})^2 - r^2) - e^{-t})dt + \sigma dW_1, \\
        dy &= (-4(y-e^{-t})((x-e^{-t})^2 + (y-e^{-t})^2 - r^2) - e^{-t})dt + \sigma dW_2,
    \end{aligned}
    \right.
\end{equation}
where $W_1, W_2$ denote independent Wiener processes. The corresponding Fokker-Planck equation is
\begin{equation}\label{eq:FP}
        \frac{\partial u}{\partial t} = \mathcal{L}u = -\nabla \cdot (fu) + \nabla^2 (Du),
\end{equation}
where $f = \begin{pmatrix}
    -4(x-e^{-t})((x-e^{-t})^2 + (y-e^{-t})^2 - r^2) - e^{-t} \\
    -4(y-e^{-t})((x-e^{-t})^2 + (y-e^{-t})^2 - r^2) - e^{-t}
\end{pmatrix}, D=\frac{\sigma^2}{2}\mathbf{I}$. The initial condition is $u(x, y, 0) = \frac{1}{K} e^{-\frac{2}{\sigma^2}((x-1)^2 + (y-1)^2 - r^2)^2}$, where $K=\iint e^{-\frac{2}{\sigma^2}((x-1)^2 + (y-1)^2 - r^2)^2} dxdy$. So the analytical solution to this equation is given by $$u(x, y, t) = \frac{1}{K} e^{-\frac{2}{\sigma^2}((x-e^{-t})^2 + (y-e^{-t})^2 - r^2)^2},$$ and the Dirichlet boundary condition are imposed by the analytical solution. In this case, we set $\sigma=0.1, r=0.5$, and the computational domain is $[0.2, 1.8]\times[0.2, 1.8]$.

\begin{figure}
    \centering
    \includegraphics[width=0.95\linewidth]{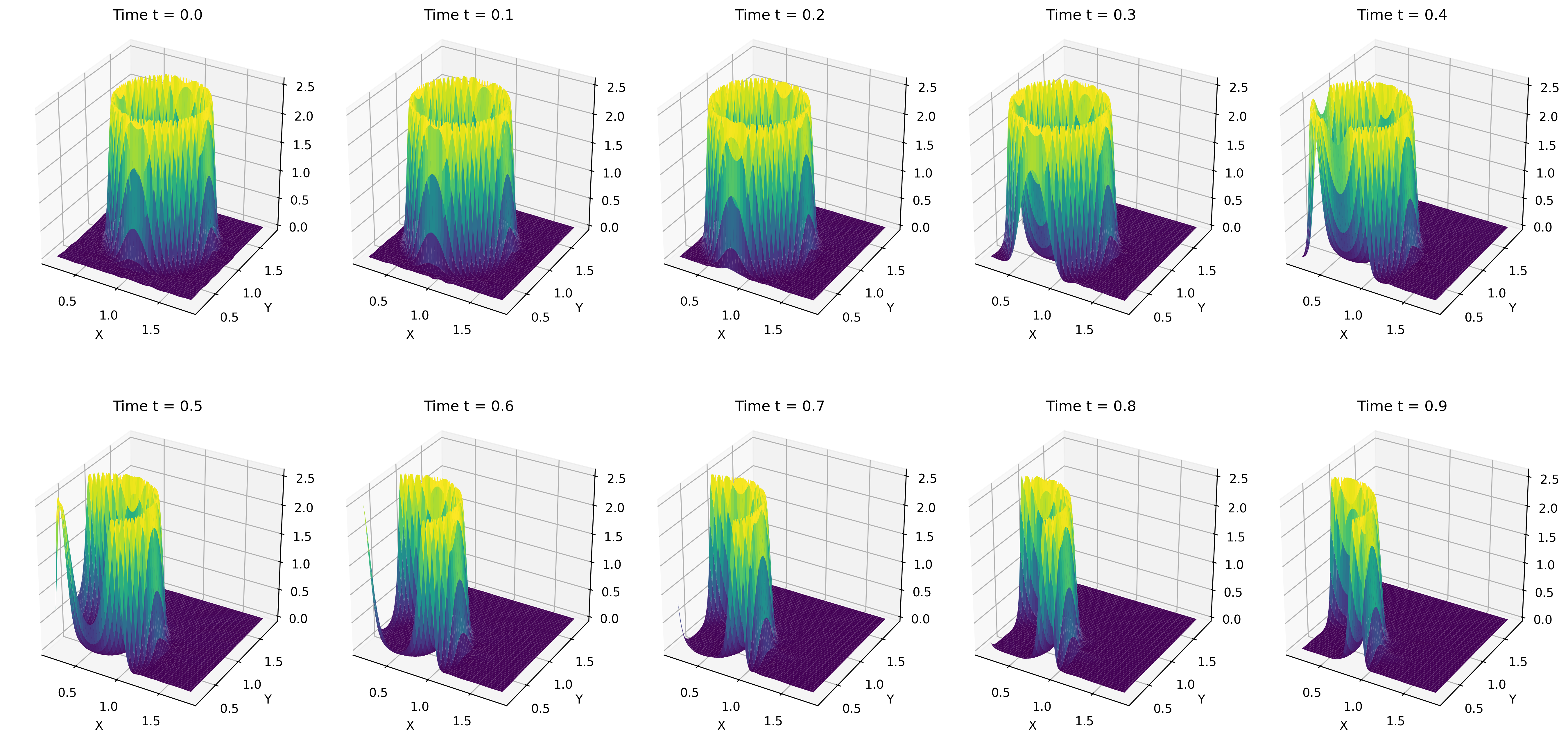}
    \caption{MSM-PINNs result for the Fokker Planck equation in Experiment~\ref{fokker planck}.}
    \label{fig:fokker planck solution}
\end{figure}

\begin{figure}
    \centering
    \includegraphics[width=0.95\linewidth]{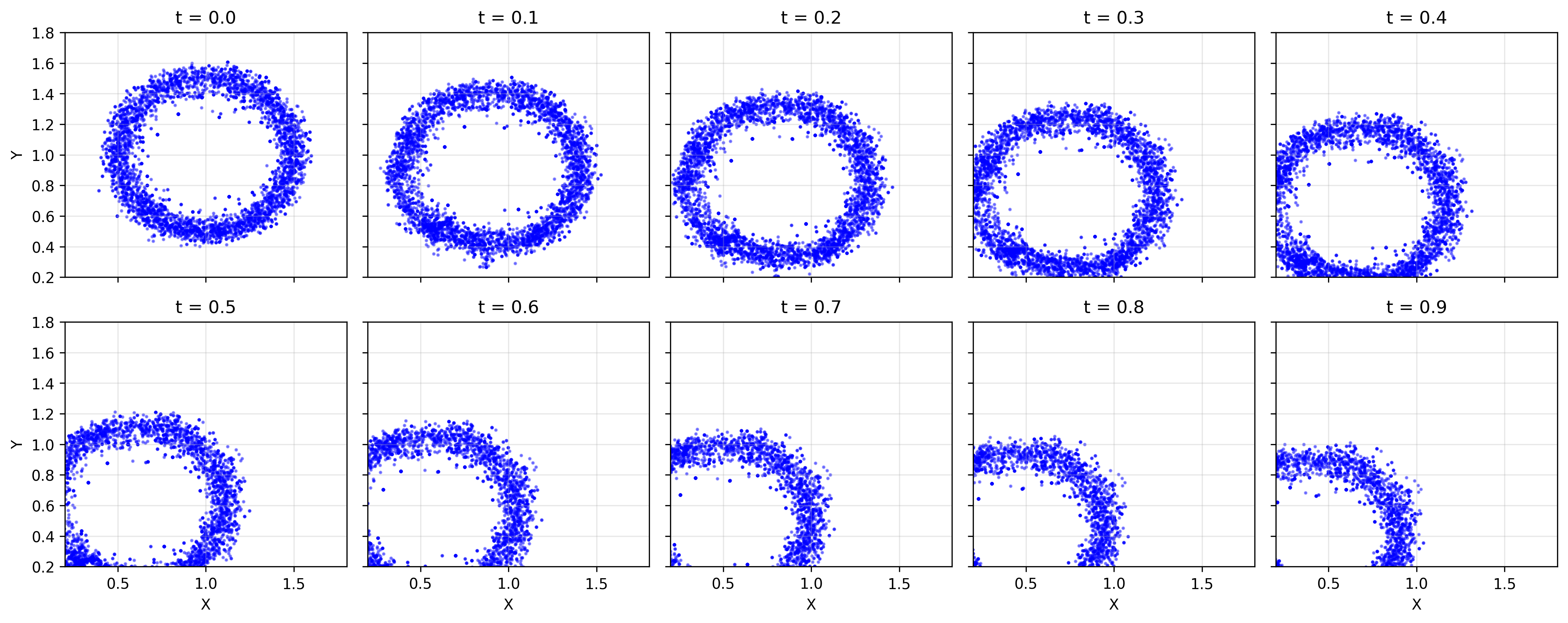}
    \caption{The trajectories of the adaptive sampling points for the Fokker Planck equation in Experiment~\ref{fokker planck}.}
    \label{fig:fokker planck sample movement}
\end{figure}

This Fokker-Planck equation governs the time evolution of the probability density for a particle undergoing Brownian motion of intensity $\sigma$ in a time-dependent ring-shaped potential centered at $(e^{-t}, e^{-t})$. The small number $\frac{\sigma^2}{2}=0.005$ in the FPE, as well as the exact solution, makes the moving distribution highly concentrated around the unit circle, which is the singularity we designed for this example around the non-trivial potential well. The deterministic drift drives the particle toward the moving potential minimum, so the density forms a ring whose center gradually shifts from $(1, 1)$ toward the origin and eventually settles into a stationary ring around $(0, 0)$. We use this example to show the movement of a relatively non-trivial distribution of singularities.

For this equation, we sample the initial training set $S_0$ proportional to the initial condition $u_0$. We set $N=1500$, $N_0=400$, $N_1=800$, $N_b=1200$ to test the performance of this example. To improve the overall coverage of the domain, we also augment the training set $S$ with additional points uniformly distributed over $\Omega \times [0, T]$, following the same strategy as in Example~\ref{Allen-Cahn}.

The solution obtained by the proposed method on the domain, together with the trajectories of the adaptive sample points, is shown in Fig.~\ref{fig:fokker planck solution} and Fig.~\ref{fig:fokker planck sample movement}. The sampling points initially lie on a ring and subsequently move together with the evolving density, remaining concentrated along the ring as it moves toward the origin. This behavior indicates that the adaptive strategy effectively tracks the motion of the highly concentrated probability mass.

\begin{figure}
    \centering
    \includegraphics[width=0.95\linewidth]{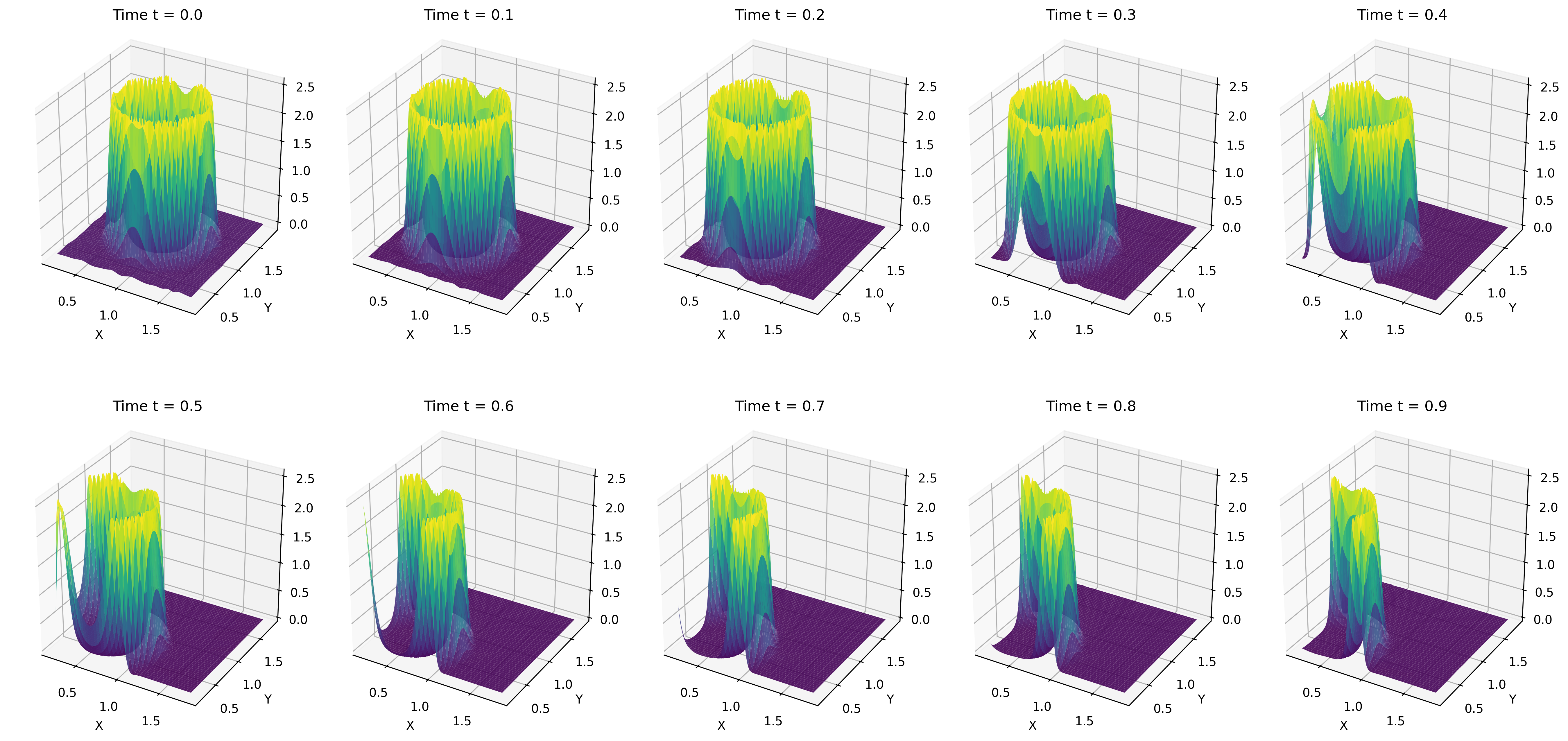}
    \caption{PINNs solution for the Fokker Planck equation in Experiment~\ref{fokker planck}.}
    \label{fig:fokker planck PINNs}
\end{figure}

Table~\ref{tab:error} summarizes the numerical errors for all methods. A representative result obtained using standard PINNs is shown in Fig.~\ref{fig:fokker planck PINNs}. In this case, the PINNs solution exhibits noticeable inaccuracies,  particularly in regions not constrained by the boundary conditions (see Fig.~\ref{fig:fokker planck PINNs}), whereas the proposed method maintains a more consistent and accurate approximation throughout the evolution.

\subsection{High-dimensional advection equation}\label{high-dim advection}

In this experiment, we consider a high-dimensional advection equation:
\begin{equation}
\left\{
\begin{array}{ll}
    \frac{\partial u}{\partial t} = -\sum_{i=1}^d \frac{\partial u}{\partial x_i}, & \text{in}\ \Omega\times[0, 1], \\
    u(\bx, 0) = e^{-\frac{1}{\alpha}\sum_{i=1}^d x_i^2}, & \text{in}\ \Omega,\\
    u(\bx, t) = e^{-\frac{1}{\alpha}\sum_{i=1}^d (x_i-t)^2}, & \text{on}\ \partial \Omega\times[0, 1],
\end{array}\right.
\end{equation}
where $d=6$ is the dimension of $\Omega$. The exact solution is $$
u(\bx, t) = e^{-\frac{1}{\alpha}\sum_{i=1}^6 (x_i-t)^2}.
$$
In this case, we set $\alpha=0.01$, $\Omega = (-0.2, 1.2)^6 \subset \mathbb{R}^6$.

\begin{figure}
    \centering
    \includegraphics[width=0.95\linewidth]{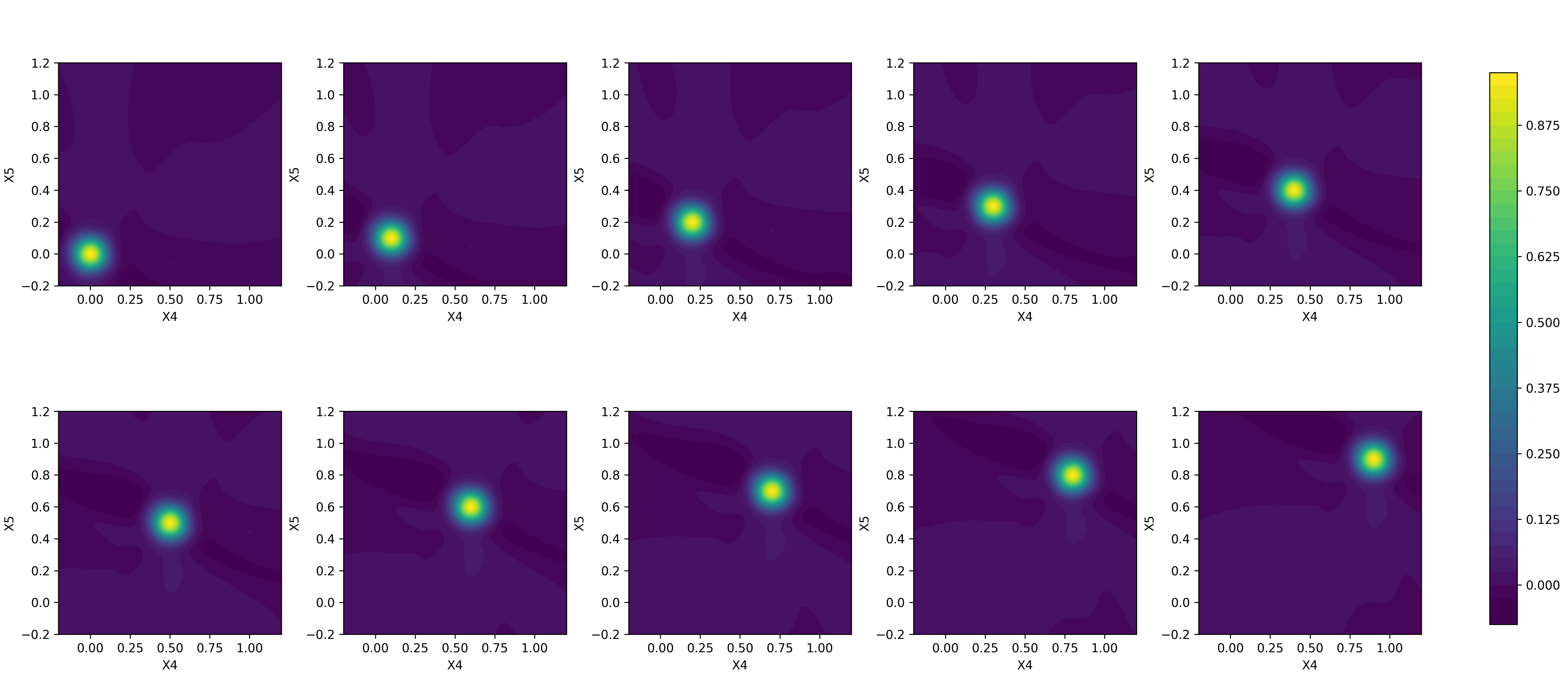}
    \caption{MSM-PINNs result for the six-dimensional advection equation in Experiment~\ref{high-dim advection}. The $i$-th subplot illustrates a two-dimensional projection of the model at time $t = 0.1(i - 1)$, with the variables $x_1$, $x_2$ $x_3$ and $x_4$ fixed at $0.1(i - 1)$.}
    \label{fig:advection}
\end{figure}

\begin{figure}
    \centering
    \includegraphics[width=0.95\linewidth]{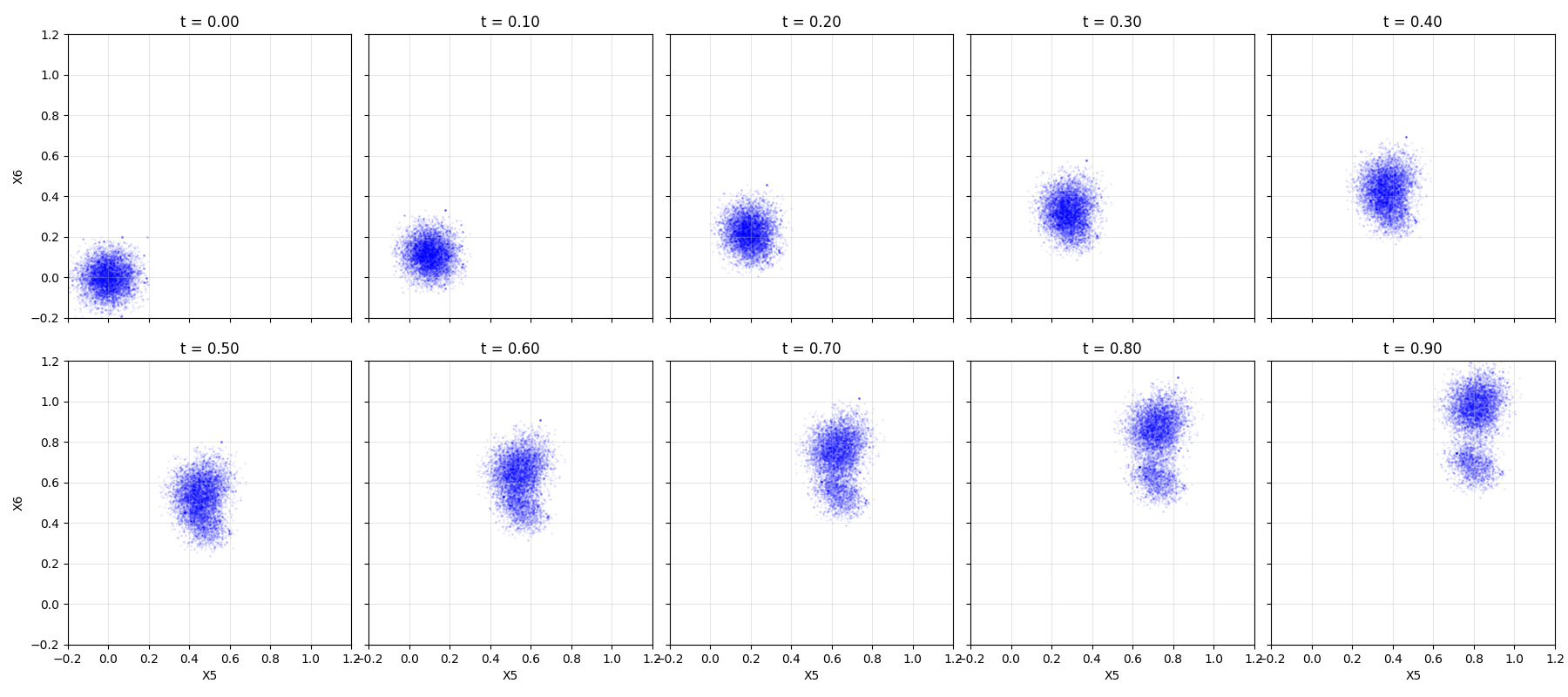}
    \caption{The trajectories of all sampling points on the last two dimensions $x_5$ and $x_6$ for the six-dimensional advection equation in Experiment~\ref{high-dim advection}.}
    \label{fig:advection sample movement}
\end{figure}

This governing equation is a $d$-dimensional linear advection equation with constant velocity vector $(1, 1, \dots, 1)$, describing rigid translation of the solution profile along the coordinate diagonal at unit speed without deformation. We set the number of sample points $N=5000, N_1=2000, N_0=2800$. For the boundary loss, we place $30$ boundary points on each coordinate hyperplane, resulting in a total of $N_b = 360$. To visualize the six‑dimensional solution, Fig.~\ref{fig:advection} shows a two-dimensional projection in the $(x_5, x_6)$ plane. The $i$-th subplot corresponds to time $t = 0.1 (i-1)$, obtained by fixing $x_1, x_2, x_3, x_4 = 0.1 (i-1)$. The projections in other pairs of dimensions are all similar. The corresponding numerical errors are reported in Table~\ref{tab:error}.

The adaptive sampling points, projected onto the last two dimensions, are shown in Fig.~\ref{fig:advection sample movement}. Most points follow the motion of the peak along the diagonal direction, which helps the network resolve the moving singularity while keeping the number of samples manageable in six dimensions. It can also be observed that in this example some sampling points do not follow the motion of the singular peak but instead drift away from it (see the last figure in Fig.~\ref{fig:advection sample movement}). This behavior stems from the use of a fixed number of iterations and a fixed number of newly added sampling points per iteration, without adapting these parameters dynamically. To be more specific, in the first iteration of adaptive sampling, the generated points have already capture the moving singularity and generate moving adaptive samples accordingly. But since $N=5000$ uniformly sampled points are still relatively sparse in this six dimensional domain (sparser than a $5^6$ mesh for FEM), the approximation away from the singularity becomes worse. %So in the second iteration of adaptive sampling, the algorithm generates new points there. The total residual of this new approximation is already much smaller than that of the first obtained with uniformly sampled collocation points. Therefore, the number of newly sampled points after the first iteration may not be as large as we fixed but should also be adaptively determined. 
As a result, in later iterations the residual near the singularity may no longer be dominant, causing the subsequently added points to shift toward other regions. In future work, we will further optimize the algorithm for such high-dimensional problems, for example, adaptively determining the number of newly sampled adaptive points in each iteration. Overall, the adaptive sampling strategy successfully captures the motion of the high-dimensional solution, whereas standard PINNs fail to approximate the problem under the same settings.

\section{Conclusion}\label{conclusion}

In this work, we proposed an adaptive sampling strategy for solving time-dependent partial differential equations, in which the distribution of sampling points at each time step is guided by the residuals $r_t(\bx)$ obtained from the previous training iteration. By concentrating subsequent training on regions with larger residuals, the network adaptively allocates computational effort where it is most needed.

Although the method requires additional training of a velocity field $\bv_{\boldsymbol{\eta}}$, which increases computational cost, it consistently achieves superior accuracy. Numerical experiments demonstrate that the approach is particularly effective for problems with singularities. In contrast to vanilla PINNs, which often rely on a large number of uniformly distributed sampling points, our method achieves higher accuracy by concentrating points in regions of strong singularity, thereby reducing the total number of required samples. This efficiency gain highlights the potential of adaptive sampling to improve both the scalability and applicability of physics-informed neural networks for complex PDEs. In particular, in a high-dimensional example, vanilla PINNs are unable to provide meaningful results with a relatively small number of collocation points, while MSM-PINNs produce reliable results with the same number of points. 

Future research will focus on optimizing the adaptive sampling framework in higher-dimensional PDE settings. In such cases, the initial sampling points used to compute the velocity field $\bv_{\boldsymbol{\eta}}$ are uniformly distributed, making it likely that none of them fall within regions of large residuals. Consequently, the adaptive sampling points may be misdirected toward less relevant regions, affecting the subsequent training process. Furthermore, when the proposed method is applied to equations in which the boundary conditions are inconsistent with the governing PDE (as in Experiment~\ref{Allen-Cahn}), the computed sampling points may drift outside the target domain $\Omega$, which can lead to numerical instability during training. At present, points that move outside the computational domain are simply discarded, which alleviates instability but reduces efficiency. Future work will therefore explore more robust handling strategies, such as projecting out-of-domain points back into feasible regions or constraining the sampling process within the domain. These improvements aim to preserve the benefits of adaptive sampling while enhancing efficiency and stability in challenging scenarios involving high dimensionality or boundary mismatches.

\section*{Acknowledgment}

This work is partly supported by the Strategic Priority Research Program of the Chinese Academy of Sciences under grant number XDA0480504 and NSFC under grant no. 12571471, 12494543, 12501598.
J. Zhai is also partially supported by the ShanghaiTech Faculty Start-Up Fund and the ShanghaiTech IMS Institutional Development Fund no. 2024X0303-902-01. X. Wan has been supported by NSF grant DMS-2513234.

\bibliography{ref}
\bibliographystyle{plain}

\clearpage
\appendix
\section{Proof of the Main Theorem}\label{appendix: proof}

To make the paper self-consistent, we provide the complete proof of the main theorem stated in the Section~\ref{dynamics} in this appendix. It is standard and can be also found in many references. All notations and assumptions are consistent with those introduced earlier.

\iffalse
\subsection{Proof of Theorem~\ref{dynamics_thm1}}

\begin{proof}
    First, by definition, we have for all $\varphi\in C^\infty_c (\Omega)$,
    $$
    \int\varphi(\bx)\,\mu_t(d\bx) = \int (\varphi\circ \bX_t(\bx))\,\nu(d\bx).
    $$
    On the other hand, we have
    $$
    \frac{\partial}{\partial t}(\varphi\circ \bX_t) = (\nabla\varphi\circ \bX_t)\cdot\frac{\partial \bX_t}{\partial t} = (\nabla\varphi\circ \bX_t)\cdot v( \bX_t(\cdot), t).
    $$
    Then for $h>0$, 
    $$\frac{\int\varphi\mu_{t+h}(d\bx)-\int\varphi\mu_{t}(d\bx)}{h} = \int\frac{\varphi\circ \bX_{t+h}(\bx) - \varphi\circ \bX_{t}(\bx)}{h}\nu(d\bx).$$
    By the assumption on $\varphi$ and the Lipschitz continuity of $T_t$, the integrand is uniformly bounded. Taking limit as $h\rightarrow 0$ and by the Dominated Convergence Theorem, $t\mapsto \int\varphi\,\mu_t(d\bx)$ is differentiable (as well as Lipschitz) for almost all $t$ on $(0,T_*)$, and
    $$
    \frac{d}{dt}\int\varphi\mu_t(d\bx) = \int(\nabla\varphi\circ \bX_t)\cdot \bv(\bX_t, t) \nu(d\bx) = \int \nabla\varphi \cdot \bv_t \,\mu_t(d\bx)
    $$
    
    Since $T_t$ is a diffeomorphism, $\mu_t=T_t\#\nu$ is also absolutely continuous with respect to the Lebesgue measure, and thus has its density function $p_t$. So
    $$
    \frac{d}{dt}\int\varphi p_t\,d\bx = \int \nabla\varphi \cdot \bv_t p_t\,d\bx = - \int \varphi \nabla\cdot(\bv_t p_t)\,d\bx
    $$
    for all $\varphi\in C^\infty_c (\Omega)$. Thus, 
    $$
    \frac{\partial p_t}{\partial t} + \nabla\cdot(p_t \bv_t) = 0, \quad 0<t<T_*
    $$
\end{proof}
\fi

\subsection{An auxiliary lemma}
We first prove the following auxiliary lemma.
\begin{lem}\label{dynamics_det}
    Let $\bv_t:\mathbb{R}^n\rightarrow\mathbb{R}^n$ be a bounded and twice differentiable velocity field of the flow mapping $\bX_t$.
    We have
    \begin{equation}
        \frac{\partial}{\partial t}\log\det\nabla_{\bx}\bX_t(\bx) = \nabla\cdot\bv(\bX_t(\bx), t).
    \end{equation}
\end{lem}

\begin{proof}
    First, 
    $$
    \bX_{t+s}(\bx) = \bX_s(\bX_t(\bx)).
    $$
    Then
    $$
    \nabla_{\bx}\bX_{t+s}(\bx) = \nabla_{\bX_t(\bx)}\bX_s(\bX_t(\bx)) \nabla_{\bx}\bX_t(\bx),
    $$
    and
    $$
    \det\nabla_{\bx}\bX_{t+s}(\bx) = \det\nabla_{\bX_t(\bx)}\bX_s(\bX_t(\bx)) \det\nabla_{\bx}\bX_t(\bx),
    $$
    So
    \begin{align*}
        \frac{\partial}{\partial t}\det\nabla_{\bx}\bX_t(\bx) & = \left[\frac{\partial}{\partial s}\det\nabla_{\bx}\bX_{t+s}(\bx)\right]_{s=0}\\
        & = \det\nabla_{\bx}\bX_t(\bx) \left[\frac{\partial}{\partial s}\det\nabla_{\bX_t(\bx)}\bX_s(\bX_t(\bx))\right]_{s=0}\\
        & = \det\nabla_{\bx}\bX_t(\bx) \left[\frac{\partial}{\partial s}\det\Big(\nabla_{\bX_t(\bx)}\bX_0(\bX_t(\bx))+s\big[\frac{\partial}{\partial s}\nabla_{\bX_t(\bx)}\bX_s(\bX_t(\bx))\big]\Big)_{s=0}\right]_{s=0}\\
        & = \det\nabla_{\bx}\bX_t(\bx) \left[\frac{\partial}{\partial s}\det\Big(\bI+s\big[\frac{\partial}{\partial s}\nabla_{\bX_t(\bx)}\bX_s(\bX_t(\bx))\big]\Big)_{s=0}\right]_{s=0}\\
        & = \det\nabla_{\bx}\bX_t(\bx)\, \text{Trace}\left( \left[\frac{\partial}{\partial s}\nabla_{\bX_t(\bx)}\bX_s(\bX_t(\bx))\right]_{s=0} \right)\\
        & = \det\nabla_{\bx}\bX_t(\bx)\, \text{Trace}\left( \nabla_{\bX_t(\bx)}\left[\frac{\partial}{\partial s}\bX_s(\bX_t(\bx))\right]_{s=0} \right)\\
        & = \det\nabla_{\bx}\bX_t(\bx)\, \text{Trace}\left(\nabla_{\bx}\bv(\bX_t(\bx), t) \right)\\
        & = \det\nabla_{\bx}\bX_t(\bx)\, \nabla_{\bx}\cdot\bv(\bX_t(\bx), t).
    \end{align*}
    Thus, 
    $$
    \frac{\partial}{\partial t}\log\det\nabla_{\bx}\bX_t(\bx) = \nabla\cdot\bv(\bX_t(\bx), t).
    $$
\end{proof}

\subsection{Proof of Theorem~\ref{dynamics_thm2}}

\begin{proof}
    The solution to equation \eqref{PDE_v_t2} exists uniquely in $\mathbb{R}^n$. If $p_t$ satisfies the transport equation \eqref{PDE_v_t2}, then it follows the ODE on the flow mapping paths
    \begin{align*}
        \frac{\partial}{\partial t}p_t(\bX_t) & =\frac{\partial}{\partial t}p_t(\bX_t)+\nabla p_t(\bX_t)\cdot\frac{\partial \bX_t}{\partial t} = \frac{\partial}{\partial t}p_t(\bX_t)+\nabla p_t(\bX_t)\cdot \bv_t(\bX_t(\cdot))\\
        & = -\nabla\cdot(p_t(\bX_t) \bv_t(\bX_t)) +\nabla p_t(\bX_t)\cdot \bv_t(\bX_t(\cdot))= -p_t(\bX_t)\nabla\cdot \bv_t(\bX_t).
    \end{align*}
    So
    $$
    p_t(\bX_t) = p_0(\bX_0)\exp\Big(-\int_0^t \nabla\cdot \bv_s(\bX_s)\,ds\Big).
    $$
    and for all test functions $\varphi\in C^\infty_c (\mathbb{R}^n)$,
    \begin{align*}
        & \int_{\mathbb{R}^n} \varphi(\bX_t(\bx)) p_0(\bX_0(\bx))\,d\bx\\
        = & \int_{\mathbb{R}^n} \varphi(\bX_t(\bx))p_t(\bX_t(\bx))\exp\Big(\int_0^t \nabla\cdot \bv_s(\bX_s(\bx))\,ds\Big)\,d\bx\\
        = & \int_{\mathbb{R}^n}\varphi(\by) p_t(\by)\,d\by,
    \end{align*}
    where we used the change of variable $\by=\bX_t(x)$ with the exponent being its Jacobian by Lemma~\ref{dynamics_det}.
    
    Conversely, we have
    $$
    \frac{\partial}{\partial t}(\varphi\circ \bX_t) = (\nabla\varphi\circ \bX_t)\cdot\frac{\partial \bX_t}{\partial t} = (\nabla\varphi\circ \bX_t)\cdot \bv( \bX_t(\cdot), t).
    $$
    Then for $h>0$, 
    $$\frac{\int\varphi\mu_{t+h}(d\bx)-\int\varphi\mu_{t}(d\bx)}{h} = \int\frac{\varphi\circ \bX_{t+h}(\bx) - \varphi\circ \bX_{t}(\bx)}{h} p_0(\bx)\,d\bx.$$
    By the assumption on $\varphi$ and the Lipschitz continuity of $\bX_t$, the integrand is uniformly bounded. Taking limit as $h\rightarrow 0$ and by the Dominated Convergence Theorem, $t\mapsto \int\varphi p_t\,d\bx$ is differentiable (as well as Lipschitz) for almost all $t$, and
    $$
    \frac{d}{dt}\int\varphi p_t\,d\bx = \int(\nabla\varphi\circ \bX_t)\cdot \bv(t, \bX_t) p_0(\bx)\,d\bx = \int \nabla\varphi \cdot p_t \bv_t\,d\bx = - \int \varphi \nabla\cdot(p_t\bv_t)\,d\bx
    $$
    Thus, 
    $$
    \frac{\partial p_t}{\partial t} + \nabla\cdot(p_t \bv_t) = 0.
    $$
\end{proof}

\section{Additional Results of the Physics-Informed Neural Networks (PINNs) Method}\label{appendix: PINNs}

This appendix provides supplementary results obtained using the Physics-Informed Neural Networks (PINNs) method for Experiments~\ref{rotation} and \ref{burgers}, intended to illustrate the instability of PINNs when solving these equations. The results shown here employ the same network architecture, number of training samples, optimizer, and number of epochs as those used for our proposed method in the main text, with modifications applied only to the weighting of the loss function components and related settings to improve the training outcomes.

\subsection{Results of the rotation equation}

\begin{figure}[h!]
    \centering
    \includegraphics[width=0.9\textwidth]{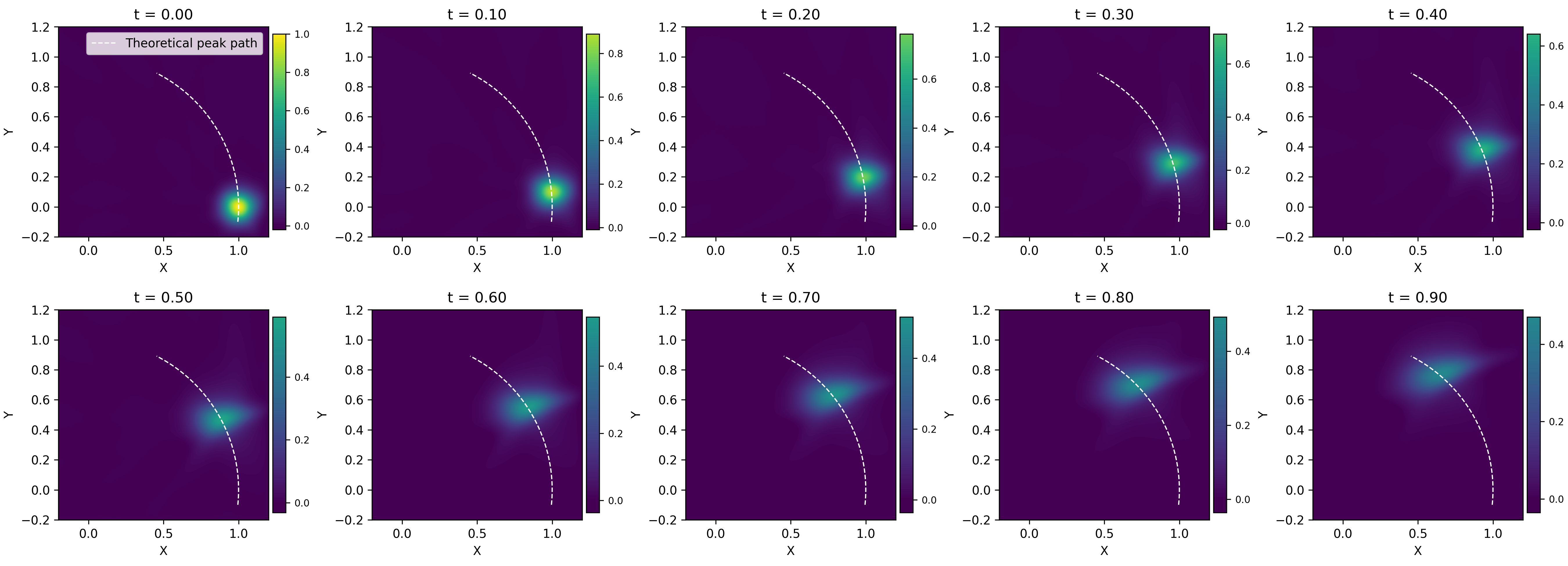}
    \caption{Numerical solution obtained by PINNs for the rotation equation in Experiment~\ref{rotation}.}
    \label{fig:exp0-1}
\end{figure}

\begin{figure}[h!]
    \centering
    \includegraphics[width=0.9\textwidth]{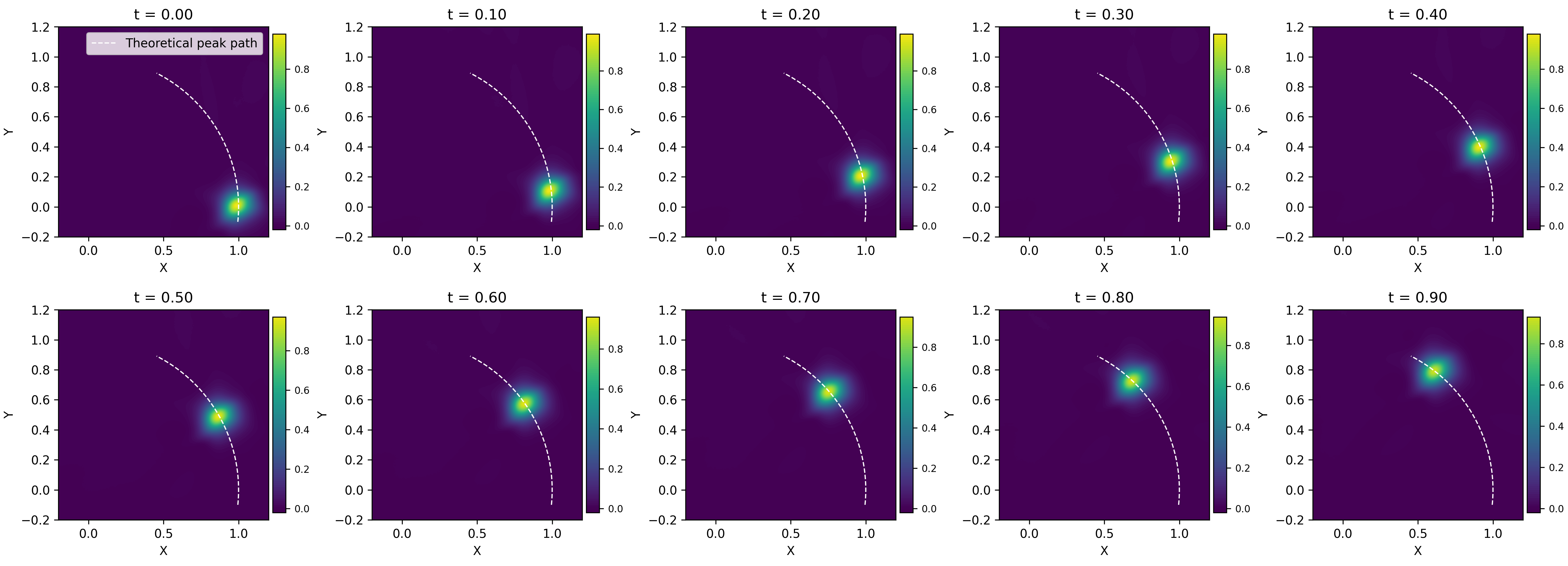}
    \caption{Numerical solution obtained by PINNs for the rotation equation in Experiment~\ref{rotation}.}
    \label{fig:exp0-2}
\end{figure}

Fig.~\ref{fig:exp0-1}--\ref{fig:exp0-2} present the results of applying the vanilla PINNs method to the rotation equation in Experiment~\ref{rotation}. When the numbers of sampling points $N, N_0$ and $N_1$ are reduced, the training process becomes highly sensitive to the choice of the balancing parameters $\beta_1, \beta_2$ in \eqref{loss_u}, and the overall behavior of PINNs becomes unstable.

Fig.~\ref{fig:exp0-1} illustrates the solution obtained with relatively large values of $\beta_1$ and $\beta_2$ in \eqref{loss_u}. Under this setting, the network fits the initial condition well; however, as time progresses, the predicted peak gradually decreases in amplitude and the originally sharp profile becomes noticeably wider. This behavior resembles numerical dissipation and indicates that the method is unable to preserve the sharp, rotating structure of the solution over time.

Fig.~\ref{fig:exp0-2} presents the result obtained with smaller loss‑term weights. In this case, the dissipation is less pronounced, but the heatmap reveals that the predicted profile deviates significantly from the expected structure. This distortion indicates that the network fails to capture the correct geometry of the solution even at the starting time.

Comparing Fig.~\ref{fig:exp0-1} and Fig.~\ref{fig:exp0-2} highlights the strong dependence of PINNs on the choice of loss-term weights: larger weights help match the initial condition but lead to severe dissipation, while smaller weights reduce dissipation but fail to capture the initial peak accurately. Such trade-offs are difficult to avoid within the PINNs framework. In contrast, our proposed method remains stable across a reasonable range of parameter choices and does not suffer from these competing failure modes.

\subsection{Results of the Burgers' equation}

\begin{figure}[h!]
    \centering
    \includegraphics[width=0.9\textwidth]{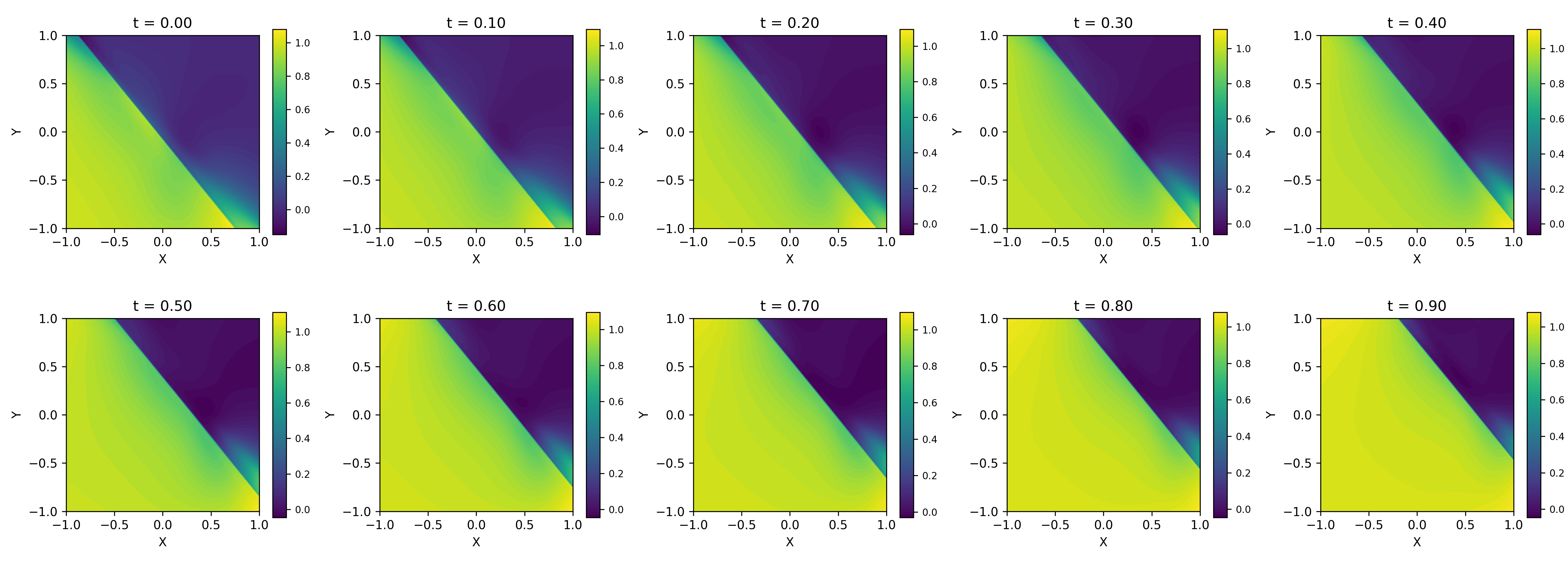}
    \caption{Numerical solution obtained by PINNs for the Burgers' equation in Experiment~\ref{burgers}.}
    \label{fig:exp1-1}
\end{figure}

\begin{figure}[h!]
    \centering
    \includegraphics[width=0.9\textwidth]{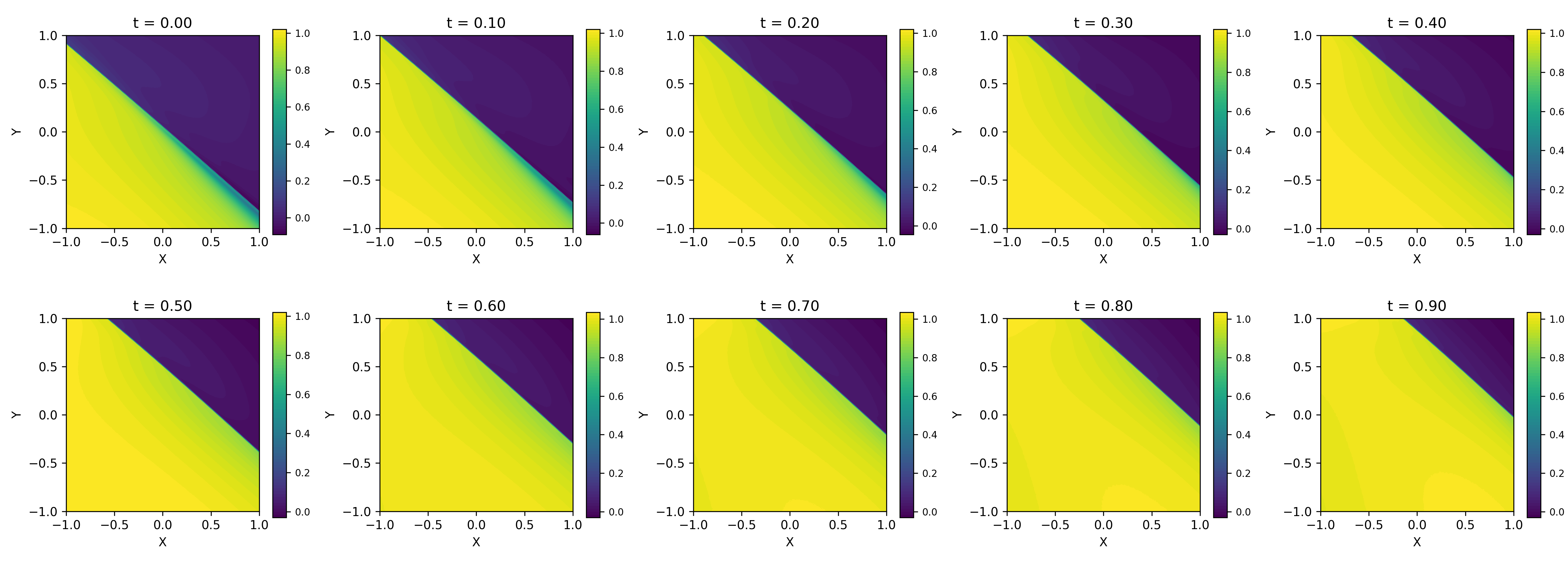}
    \caption{Numerical solution obtained by PINNs for the Burgers' equation in Experiment~\ref{burgers}.}
    \label{fig:exp1-2}
\end{figure}

Fig.~\ref{fig:exp1-1}--\ref{fig:exp1-2} show two results obtained by applying the vanilla PINNs to the Burgers' equation in Experiment~\ref{burgers}. When the number of sampling points is reduced, the performance of PINNs becomes unstable and the predicted solutions vary noticeably across different runs.

Both Fig.~\ref{fig:exp1-1} and Fig.~\ref{fig:exp1-2} exhibit a common failure mode: the predicted interface is clearly shifted relative to the analytical solution, and the network fails to capture the sharp gradient transition characteristic of the Burgers' equation. Although the overall shape of the solution seems reasonable, the displacement of the interface indicates that vanilla PINNs struggle to resolve steep features accurately under this setting.

These examples highlight the difficulty of applying vanilla PINNs to problems with sharp interfaces when only a limited number of sampling points is available. In contrast, our method remains stable within a reasonable range of parameter choices and does not exhibit such sensitivity.

\end{document}